\documentclass[3p,preprint]{elsarticle}

\usepackage{amsmath}
\usepackage{amsfonts}
\usepackage{amssymb,latexsym}
\usepackage{amsthm}
\usepackage{graphicx}
\usepackage{ulem}
\usepackage{color}
\usepackage{newlfont} 
\usepackage{subfigure}
\usepackage{verbatim}
\usepackage{rotating}
\usepackage{multirow}
\usepackage{units}
\usepackage{bm}
\usepackage[english]{babel}
\usepackage[utf8]{inputenc}
\usepackage{algorithm}
\usepackage[noend]{algpseudocode} 
\usepackage{booktabs}
\usepackage{caption}
\usepackage{hyperref}
\usepackage{subfigure}
\usepackage[export]{adjustbox}

\newtheorem{theorem}{Theorem}[section]

\newtheorem{Remark}{Remark}[section]

\newtheorem{Proposition}{Proposition}[section]

\newcommand{\lambdamax}{\lambda_{\mbox{\tiny{max}}}}
\newcommand{\lambdamin}{\lambda_{\mbox{\tiny{min}}}}
\newcommand{\condmin}{{{\mbox{{cond}}}_2({\mbox{G}}}_{n_{min}})}
\newcommand{\condmax}{{{\mbox{{cond}}}_2({\mbox{G}}}_{n_{max}})}

\journal{...}

\begin{document}

\begin{frontmatter}

\title{{\color{black}On the role of weak Marcinkiewicz-Zygmund constants in\\ polynomial approximation by orthogonal bases}}

\author[addressHK]{Congpei An}
 
\author[address-PD]{Alvise Sommariva}
\author[address-PD]{Marco Vianello\corref{corrauthor}}
\ead{marcov@math.unipd.it}

\cortext[corrauthor]{Corresponding author}

\address[addressHK]{School of Mathematics and Statistics, Guizhou University, Guiyang 550025, Guizhou, China}
\address[address-PD]{Department of Mathematics "Tullio Levi-Civita", University of Padova, Italy}

\begin{abstract}
{\color{black}
We compute numerically the $L^2$ Marcinkiewicz-Zygmund constants of cubature rules, with a special attention to their role in polynomial approximation by orthogonal bases. 
We test some relevant rules on domains such as the interval, the square, the disk, the triangle, the cube and the sphere.
The approximation power of the corresponding least squares (LS) projection is compared with standard hyperinterpolation and its recently proposed ``exactness-relaxed'' version. The Matlab codes used for these tests are available in open-source form.}

\end{abstract}


\end{frontmatter}

\section{Introduction}
In this paper we compute numerically the $L^2$ Marcinkiewicz-Zygmund constants of cubature rules $S$ on multivariate domains $\Omega$, with a special attention to their role in polynomial approximation by orthogonal bases.
We recall that a cubature rule with positive weights $w_i>0$
\begin{equation} \label{S}
S(f)=\sum_{i=1}^M{w_i\,f(\mathbf{x}_i)}\approx I(f)=\int_\Omega{f(\mathbf{x})\,d\mu}
\end{equation}
has the MZ (Marcinkiewicz-Zygmund) property over $\Omega$ {\color{black}for a measure $\mu$}, if there exist constants $A>0$ and $B > 0$ such that 
\begin{equation} \label{MZD}
A \|p\|_{L^2}^2 \leq S(p^2)  \leq B \|p\|_{L^2}^2, \quad  \forall p \in {\mathbb{P}}_n\;,\;\;\|p\|_{L^2}^2=I(p^2)\;.
\end{equation}
In general is of interest to determine the larger $A=A(n)>0$ and smaller $B=B(n)>0$ for which (\ref{MZD}) is verified. We observe that (\ref{MZD}) is a MZ inequality in {\em weak form}, since the constants $A$ and $B$ are independent of $p\in \mathbb{P}_n$ but dependent on $n$. Indeed, starting from the original Marcinkiewicz-Zygmund paper in 1937 \cite{MZ37}, the study of conditions ensuring that such constants are bounded from below and above, respectively, 
i.e. (\ref{MZD}) holds with constants independent of $n$, has been object of a specific literature and is still an active research topic; with no pretence of exhaustivity, we may quote e.g. \cite{DKP25,FM10,L98} with the references therein. 

In the present paper we focalize on the MZ inequality in its weak form (\ref{MZD}) and on the computation of the MZ constants. Indeed, 
given a $\mu$-orthonormal basis $\{\phi_j\}$ of $\mathbb{P}_n$, as claimed in \cite{ANWU22} if 
\begin{equation} \label{eta}
\eta=\eta(n)=\max\{|A-1|,|B-1|\}<1 
\end{equation}
and it is not too close to 1, that is if $B<2$ and $A>0$ is not too small, {\color{black}by the accessory MZ inequality
\begin{equation} \label{etabound}
|S(p^2)-I(p^2)|\leq \eta\,I(p^2)\;,\;\;\forall p\in \mathbb{P}_n\;,
\end{equation}
it follows that the ``exactness-relaxed'' hyperinterpolation operator
\begin{equation} \label{hyper} 
\mathcal{H}_nf=\sum_{j=1}^{d_n}{S(f\phi_j)\,\phi_j}\;,\;\;d_n=dim(\mathbb{P}_n)\;,
\end{equation}
can be used for the approximation of $f\in C(\Omega)$, even though the rule does not integrate exactly in 
$\mathbb{P}_{2n}$ and thus $\mathcal{H}_n$ is not a projection.} Namely, in case $S$ is an algebraic cubature rule exact in $\mathbb{P}_{n+k}$ with $k\leq n$, in \cite{ANWU22} it is proved that 
\begin{equation} \label{n+k-err}
\|f-\mathcal{H}_nf\|_{L^2}\leq \left(1+\frac{1}{\sqrt{1-\eta}}\right)\,\sqrt{\mu(\Omega)}\;E_k(f)\;.
\end{equation}
Notice that for $k=n$ we have $\eta=0$ and we recover standard Sloan's hyperinterpolation \cite{Slo95}. The method is also called ``unfettered'' hyperinterpolation in \cite{ANWU24}.
On the other hand, {\color{black}when $S$ is a QMC (Quasi-MonteCarlo) rule} then, as recently proved in \cite{AKG25}, {\color{black}if $0<\eta<1$}
\begin{equation} \label{QMC-err}
\|f-\mathcal{H}_nf\|_{L^2}\leq \left(1+\sqrt{1+\eta}\right)\,\sqrt{\mu(\Omega)}\,E_n(f)+\|p^\ast_n-\mathcal{H}_np^\ast_n\|_{L^2}\;,
\end{equation}
where $p^\ast_n$ is the best uniform approximation to $f$ in $\mathbb{P}_n$, i.e., $\|f-p^\ast_n\|_\infty=E_n(f)$. The second summand on the r.h.s. of (\ref{QMC-err}), which is not zero because $\mathcal{H}_n$ is not a projection operator in this case, is termed in \cite{AKG25} the ``aliasing term'' of unfettered hyperinterpolation. {\color{black}For these recent developments of hyperinterpolation theory we refer the reader also to the survey paper \cite{ARW25}.}

On the other hand, as observed for example in \cite{FM10,G20} and proved in detail for completeness in the next section, $A$ and $B$ are nothing but the smallest and largest eingenvalue of the Gramian of the discrete measure corresponding to the cubature rule. Thus, 
under the conditions above on $A$ and $B$, the Gramian is well-conditioned and we can also compute with high accuracy the discrete orthogonal projection on $\mathbb{P}_n$
\begin{equation} \label{orth} 
\mathcal{G}_nf=\sum_{j=1}^{d_n}{c_j\,\phi_j}\;,\;\{c_j\}=G^{-1}\{S(f\phi_j)\}\;,\;\;G=(S(\phi_i\phi_j))_{1\leq i,j\leq d_n}\;.
\end{equation}

{\color{black} Concerning the Least-Squares approximation error, we can state and prove the following result. 
\begin{Proposition}
Let $\mathcal{G}_nf$ be the weighted Least-Squares polynomial defined by (\ref{S}) and (\ref{orth}), where $f\in C(\Omega)$ and the cubature rule satisfies the MZ property (\ref{MZD})-(\ref{eta}). Then the following bound holds 
\begin{equation} \label{alt}
\|f-\mathcal{G}_nf\|_{L^2} 
\leq \left(1+\sqrt{\mbox{cond}_2(G)}\,\right) \sqrt{\mu(\Omega)} \,E_n(f)
\leq \left(1+\sqrt{\frac{1+\eta}{1-\eta}}\right)\,\sqrt{\mu(\Omega)}\;E_n(f)\;,
\end{equation}
which, in case the cubature rule is exact on the constants, can be refined to
\begin{equation} \label{err-Gn}
\|f-\mathcal{G}_nf\|_{L^2}\leq 
\left(1+\frac{1}{\sqrt{A}}\,\right)\,\sqrt{\mu(\Omega)}\,E_n(f)\leq 
\left(1+\frac{1}{\sqrt{1-\eta}}\right)\,\sqrt{\mu(\Omega)}\;E_n(f)\;.
\end{equation}
\end{Proposition}

\noindent{\bf Proof.} In view of the MZ inequality (\ref{MZD}) and the fact that $\mathcal{G}_nf$ is an orthogonal projection, we can write the chain of inequalities
\begin{align*}
\| f - \mathcal{G}_n f \|_{L^2} 
&\le \| f - p_n^* \|_{L^2} + \| p_n^* - \mathcal{G}_n p_n^* \|_{L^2} + \| \mathcal{G}_n (p_n^* - f) \|_{L^2} \\
&= \| f - p_n^* \|_{L^2} + \| \mathcal{G}_n (p_n^* - f) \|_{L^2} \\
&\le \| f - p_n^* \|_{L^2} + \frac{1}{\sqrt{A}} \sqrt{S \bigl( (\mathcal{G}_n(p_n^* - f))^2 \bigr)} \\
&= \| f - p_n^* \|_{L^2} + \frac{1}{\sqrt{A}} \| \mathcal{G}_n (p_n^* - f) \|_{\ell^2_w} \\
&\le \| f - p_n^* \|_{L^2} + \frac{1}{\sqrt{A}} \| p_n^* - f \|_{\ell^2_w}
\end{align*}
\begin{equation} \label{chain}
\le \biggl( \sqrt{\mu(\Omega)} + \sqrt{\frac{\sum w_i}{A}} \; \biggr) \| f - p_n^* \|_\infty
= \biggl( \sqrt{\mu(\Omega)} + \sqrt{\frac{\sum w_i}{A}} \; \biggr) E_n(f).
\end{equation}
Now, if the cubature rule is exact at least on the constants (as any algebraic rule or also a QMC rule \cite{DP10}), we have $\sum{w_i}=\mu(\Omega)$ and thus (\ref{err-Gn}) holds
$$
\|f-\mathcal{G}_nf\|_{L^2}\leq 
\left(1+\frac{1}{\sqrt{A}}\,\right)\,\sqrt{\mu(\Omega)}\,E_n(f)\leq 
\left(1+\frac{1}{\sqrt{1-\eta}}\right)\,\sqrt{\mu(\Omega)}\;E_n(f)\;.
$$

Alternatively, by the MZ inequality, since the weights are positive,
$$
\sum w_i=|S(1)| \leq B \|1\|_{L^2}^2 =B \int_\Omega 1 \, d\mu=B \mu(\Omega)
$$
and thus, from $\mbox{cond}_2(G)=B/A$ and the last row of (\ref{chain}), we get (\ref{alt}) 
$$
\|f-\mathcal{G}_nf\|_{L^2} 
\leq \left(1+\sqrt{\mbox{cond}_2(G)}\,\right) \sqrt{\mu(\Omega)} \,E_n(f)
\leq \left(1+\sqrt{\frac{1+\eta}{1-\eta}}\right)\,\sqrt{\mu(\Omega)}\;E_n(f)\;,
$$
which holds even if the cubature rule $S$ is not exact on the constants. \hspace{0.5cm}$\square$
\vskip0.5cm

\begin{Remark}
Notice that if $S$ is an algebraic cubature rule exact in $\mathbb{P}_{n+k}$ with $k<n$, then (\ref{err-Gn}) is a better estimate than  
(\ref{n+k-err}). Moreover, for a general cubature rule satisfying (\ref{MZD})-(\ref{eta}), both (\ref{alt}) and (\ref{err-Gn}) do not contain any aliasing term, differently from the unfettered hyperinterpolation (\ref{QMC-err}), because $\mathcal{G}_nf$ is a projection.
\end{Remark}

The paper is organized as follows.} In the next section, we start showing how to compute the $L^2$ Marcinkiewicz-Zygmund costants $A$, $B$ for a fixed cubature rule $S$, as extremal eigenvalues of the Gramian corresponding to the rule. Next, we prove that univariate $k$-points Gaussian rules do not satisfy (\ref{eta}) for $n\geq k$.

In the numerical section we determine the values of MZ constants for some relevant algebraic rules in domains as the unit interval, square, disk, triangle and sphere. {\color{black}In particular we try to answer numerically the following question: 
\begin{itemize}
\item {\em a rule with Algebraic Degree of Exactness $ADE=m$ can have $\eta < 1$ for some ${\color{black}n>\lfloor\frac{m}{2}\rfloor}$} ?
\end{itemize}
or in other words, by a dual point of view, 
\begin{itemize}
\item {\em for a fixed degree $n$, a rule with $ADE=m<2n$ can have $\eta < 1$} ?
\end{itemize}
}

This cannot never hold for Gaussian rules, as proved theoretically in Section 2. 
For the analyzed degrees this is instead the case of Clenshaw-Curtis rule for the interval, of Padua points based rule in the square as well as of some spherical and symmetric spherical designs proposed in \cite{RWH}, in which $\eta$ is not too close to $1$, suggesting their use in (\ref{hyper}) and (\ref{orth}). The feature comes out also with some instances of near minimal rules on the disk and the triangle. In these cases however $\eta \approx 1$, showing that their usage in (\ref{hyper}) and (\ref{orth}) can be considered but not suggested.

{\color{black}
Moreover, for any convergent cubature rule we have that $A,B\to 1$ as $M\to \infty$ for fixed $n$, since the Gramian converges to the identity matrix, so $\eta\to 0$ and thus is bounded away from 1 for $M=M(n)$ sufficiently large. This is for example the case of QMC rules, that we will test numerically in the square and the cube for a range of degrees.} 

All the Matlab routines used in this work and the demos are available as open-source codes at {\cite{S25}}. {\color{black}We stress that the aim of the present paper is not to study theoretically or to compute numerically (if feasible) ``strong'' (i.e., independent of $n$) MZ constants, topics for which we refer the reader to the specialized literature, partially quoted in the bibliography. We try instead to give some computational tools, 
including the related codes, that allow to check the approximation power of polynomial approximation via orthogonal bases at a given degree $n$, via ``weak'' (i.e. dependent on $n$) MZ constants of widely adopted cubature rules.}

\section{On the computation of weak Marcinkiewicz-Zygmund constants} 

The aim of this section is to show how to compute numerically {\color{black}{the largest and smallest nonnegative constants $A=A(n)$ and $B=B(n)$}} such that 
\begin{equation}\label{MZDA}
A \|p\|_{L^2}^2 \leq |S(p^2)|  \leq B \|p\|_{L^2}^2, \quad  \forall p \in {\mathbb{P}}_n\;.
\end{equation}
{\color{black}To  this purpose, for convenience we can state and prove the following fundamental result, for which we refer also to \cite{FM10,G20}. 

\begin{Proposition}[Eigenvalue Characterization of Weak Marcinkiewicz--Zygmund Constants]
Let \(S\) be a positive-weight cubature rule, and let \(\{\phi_j\}_{j=1}^{d_n}\) be a 
\(\mu\)-orthonormal basis of the polynomial space \(\mathbb{P}_n\).  
Denote by \(G\) the corresponding Gram matrix with entries  

\begin{equation} \label{eq:Gram}
G_{j,k} = S(\phi_j \phi_k), \quad 1 \le j,k \le d_n.
\end{equation}
Then the optimal constants \(A = A(n)\) and \(B = B(n)\) satisfying the weak Marcinkiewicz--Zygmund inequality  

\begin{equation} \label{eq:MZ-ineq}
A \|p\|_{L^2}^2 \le S(p^2) \le B \|p\|_{L^2}^2, \qquad \forall p \in \mathbb{P}_n\;,
\end{equation}
are precisely the smallest and largest eigenvalues of the Gram matrix \(G\):
\begin{equation} \label{eq:eigenvalues}
A = \lambda_{\min}(G), \qquad B = \lambda_{\max}(G).
\end{equation}

Consequently, the constant $\eta$ in (\ref{etabound}) is 
\begin{equation} \label{eq:eta}
\eta = \max\bigl\{ |1-A|, |1-B| \bigr\}\;.
\end{equation}
\end{Proposition}
}
\noindent{\bf Proof.}
{\color{black}Since inequality (\ref{eq:MZ-ineq}) is trivially satisfied for the null polynomial with any $A,B>0$}, one can restrict the problem above to the case in which $\|p\|_{L^2}=1$, 
that is to find 
\begin{equation}\label{MZAB}
A=\min_{p \in {\mathbb{P}}_n, \|p\|_{L^2}=1} |S(p^2)|\;,\;\;
B=\max_{p \in {\mathbb{P}}_n, \|p\|_{L^2}=1} |S(p^2)|\;. 
\end{equation}
Assume that $p=\sum_{j=1}^{d_n} c_j \phi_j$ is such that $\|p\|_{L^2}^2=I(p)=1$, that is $\sum_{j=1}^{d_n} c^2_j=1$.
For $\mathbf{c}=\{c_j\}$, we have that
\begin{eqnarray}\label{4}
S(p^2)&=&\sum_{i=1}^{M} w_i p^2({\mathbf{x}}_i) = \sum_{i=1}^{M} w_i \left(\sum_{j=1}^{d_n} c_j \phi_j ({\mathbf{x}}_i)\right)^2 \nonumber \\
&=&\sum_{i=1}^{M} w_i \left(\sum_{j=1}^{d_n} c^2_j \phi^2_j ({\mathbf{x}}_i)  +2   \sum_{j,k=1, j< k}^{d_n} c_j c_k \phi_j ({\mathbf{x}}_i) \phi_k ({\mathbf{x}}_i)\right) \nonumber \\
&=& \sum_{j=1}^{d_n} c^2_j \sum_{i=1}^{M} w_i \phi^2_j ({\mathbf{x}}_i) +2   \sum_{j,k=1, j< k}^{d_n} c_j c_k \sum_{i=1}^{M} w_i \phi_j ({\mathbf{x}}_i) \phi_k ({\mathbf{x}}_i)  \nonumber \\
&=&\sum_{j=1}^{d_n} c^2_j S (\phi^2_j) + 2   \sum_{j,k=1, j< k}^{d_n} c_j c_k S (\phi_j\phi_k)  \nonumber \\ &=&\sum_{j=1}^{d_n} c^2_j G_{j,j} + 2   \sum_{j,k=1, j< k}^{d_n} c_j c_k G_{j,k}  \nonumber \\
&=& {{\mathbf{c}}^T G {\mathbf{c}}}\;.
\end{eqnarray}

Thus the largest $A$ and the smallest $B$ satifying (\ref{MZAB}) correspond to determine the solution of the optimization problems
\begin{equation}\label{AB}
\min_{\mathbf{c}^T \mathbf{c}=1}{{\mathbf{c}}^T G {\mathbf{c}}}\;,\;\;\;
\max_{\mathbf{c}^T \mathbf{c}=1}{{\mathbf{c}}^T G {\mathbf{c}}}\;.
\end{equation}
In other words, we need to minimize/maximize the quadratic form associated with the Gramian $G$ on the unit-sphere. By the {\it{min-max}} theorem, 
\begin{itemize}
\item $A$ is the minimum eigenvalue $\lambdamin(G)$ of the positive semi-definite matrix $G$,
\item $B$ is the minimum eigenvalue $\lambdamax(G)$ of the positive semi-definite matrix $G$.
\end{itemize}
We observe that $A$ may be $0$, since $G$ is only positive semi-definite. In this case the rule does not possess a Marcinkiewicz-Zygmund property, since it is required that $A>0$.
Moreover, the polynomials $p_A$ and $p_B$, whose coefficients in the orthonormal basis are the components of the unitary eigenvector relative to the extremal eigenvalues of $G$, are such that 
\begin{equation}
A= S(p_A^2)\;,\;\;B=S(p_B^2)\;.
\end{equation}

We now focus our attention in determining the smallest $\eta \geq 0$ such that
\begin{equation}\label{MZ}
|I(p^2)-S(p^2)| \leq \eta |I(p^2)|
\end{equation}
for any $p \in {\mathbb{P}}_{n}$. Again, it is easy to see that we can restrict our attention to $p \in \mathbb{P}_n$ that are unitary in norm $2$, so obtaining
\begin{equation}\label{MZmax}
\eta=\max_{p \in {\mathbb{P}}_n, \|p\|_{L^2}=1} |I(p^2)-S(p^2)|\;.
\end{equation}

Now consider $p=\sum_{j=1}^{d_n} c_j \phi_j$ and suppose $I(p^2)=1$, that is $\sum_{j=1}^{d_n} c^2_j=1$. From (\ref{4}), 
\begin{eqnarray}\label{in}
|I(p^2)-S(p^2)|&=&|1-S(p^2)|  \nonumber \\ 
&=&\left|1-\left(\sum_{j=1}^{d_n} c^2_j S (\phi^2_j) + 2   \sum_{j,k=1, j< k}^{d_n} c_j c_k S (\phi_j\phi_k)\right)\right|  \nonumber \\
&=&\left|1-\left(\sum_{j=1}^{d_n} c^2_j G_{j,j} + 2   \sum_{j,k=1, j< k}^{d_n} c_j c_k G_{j,k}\right)\right|  \nonumber \\
&=& |1 - {\mathbf{c}}^T G {\mathbf{c}}|=|{\mathbf{c}}^T  {\mathbf{c}} - {\mathbf{c}}^T G {\mathbf{c}}| \nonumber \\
&=& |{\mathbf{c}}^T (Id-G) {\mathbf{c}}|\;,
\end{eqnarray}
and to obtain $\eta$ we have to maximize the r.h.s. of ({\ref{in}}) for all vectors ${\mathbf{c}}$ such that $\|{\mathbf{c}}\|_{L^2}=1$, that corresponds to compute the spectral radius of the error matrix $E=Id-G$, where $Id$ is the identity matrix. 
In other words
\begin{equation}\label{etaAB}
\eta={\color{black}{\|E\|_2=}} \max\{|1-\lambdamin (G)|,|1-\lambdamax (G)|\}=\max\{|1-A|,|1-B|\}\;.\hspace{0.5cm}\square
\end{equation}

\vskip0.5cm

\begin{Remark}
Observe that the result above  does not assume that the formula $S$ has a certain degree of exactness  and can be applied to general measures $\mu$, even discrete ones.
\end{Remark}

\begin{Remark}
{\color{black}Assume that a sequence of formulas $S_m$, $m=1,2,\ldots$ is convergent on the polynomials, which as known by the multivariate generalization of Polya-Steklov theorem is a necessary and sufficient condition for convergence in $C(\Omega)$. Then the corresponding sequence of Gram matrices $G_m$ converges to the identity matrix, and their eigenvalues coinverge to 1. Consequently, fixed arbitrarily $\eta\in (0,1)$ there exists $m^* = m^*(n)\in {\mathbb{N}}$ such that $S_{m^*}$ satisfies the property ({\ref{MZ}}) for $p \in {\mathbb{P}}_n$.
An example is that of QMC rules based on an increasing number of low-discrepancy points on $\Omega$, since they are convergent on any fixed polynomial and in particular on the polynomials $\phi_j\phi_k\in \mathbb{P}_{2n}$, $1\leq j,k\leq d_n$.}
\end{Remark}

\begin{Remark}
Let $\hat{c}\in {\mathbb R}^{d_n}$ be a unitary eigenvector of $E=Id-G$ w.r.t. to the eigenvalue $\lambdamax (E)$ of largest magnitude of $E$. From the interpretation of the problem in linear algebra terms, we have that the polynomial $\hat{p} \in {\mathbb{P}}_n$ of 2-norm equal to 1, that achieves the worst case error in (\ref{MZ}) is $\hat{p}=\sum_{j=1}^{d_n} \hat{c}_j \phi_j$.
\end{Remark}

\subsection{A negative result about Gaussian rules}

{\color{black}
As recalled in the Introduction, it is interesting to check whether a quadrature rule with Algebraic Degree of Exactness $ADE=m$ could have $\eta < 1$ for some $n>\lfloor\frac{m}{2}\rfloor$; cf. (\ref{eta})-(\ref{n+k-err}). Indeed, in such a case it would be possible to apply for example exactness-relaxed hyperinterpolation as proposed in \cite{ANWU22}.}

We investigate here the case of {\color{black}$k$-points Gaussian rules $S_k$, that have ADE $m=2k-1$} w.r.t. a weight function $w$ on the interval $(a,b)$, {\color{black}possibly unbounded}.
It is obvious that if $n \in {\mathbb N}$ and $n < k$ then $2n < 2k -1$, easily implying in view of the ADE of the rule that $\eta = 0$. 

Take $n=k$. Let $p=\sum_{i=0}^{n} c_i \psi_i$ where $\{\psi_i\}_{i=0,\ldots,n}$ is the triangular family of orthonormal polynomials w.r.t. the weight function $w$, that is ${\mbox{deg}}(\psi_i)=i$ and $\int_a^b \psi_i(x) \psi_j(x) w(x) dx =\delta_{i,j}$ for $i,j=0,\ldots,n$.
Denoting by $\{\psi^*_i\}_{i=0,\ldots,n}$  the triangular family of monic orthogonal polynomials w.r.t. the weight function $w$ in $(a,b)$ we observe that if $\psi_{n}(x)=\alpha_n x^n +\ldots + \alpha_0$ then $\psi_{n}(x)=\alpha_n \psi^*_{n}(x)$ since $\psi^*_{n}$ is monic. In particular $\alpha_n=\frac{1}{\|\psi^*_{n}\|_{L^2_w}}$.
We recall that if $f \in C^{2n}(a,b)$ then 
\begin{equation}\label{EG1}
\int_a^b f(x) w(x) dx -S_n(f) =\frac{D^{(2n)}f(\xi)}{2n!} \|\psi^*_n\|^2_{L^2_w}, \quad \xi \in (a,b),
\end{equation}
in which $\|f\|_{L^2_w}=\int_a^b f^2(x) w(x) dx$ for $f \in L_{L^2_w}$ (see, e.g., {\cite[p.159]{SB10}}).

Now 
\begin{itemize}
\item suppose that $\|p\|_{L^2_w}=1$ or equivalently $\sum_{i=0}^{n} c^2_i=1$;
\item set $q_{n-1}=\sum_{i=0}^{n-1} c_i \psi_i \in {\mathbb{P}}_{n-1}$ so that $p=c_{n} \psi_{n} + q_{n-1}$.
\end{itemize}
In view of ({\ref{EG1}}), for some $\xi \in (a,b)$, we get 
\begin{eqnarray}\label{EG2}
\int_a^b p^2(x) w(x) dx -S_n( p^2 ) &=&\frac{ D^{(2n)}p^2(\xi)}{2n!} \|\psi^*_n\|^2_{L^2_w} = \frac{ D^{(2n)}(c_{n} \psi_{n} + q_{n-1})^2 (\xi)}{2n!} \|\psi^*_n\|^2_{L^2_w} \nonumber \\
&=& c^2_{n}  \frac{ D^{(2n)} \psi^2_{n}(\xi)}{2n!} \|\psi^*_n\|^2_{L^2_w} = c^2_{n}  \frac{ 2n! \alpha^2_n }{2n!} {\frac{1}{\alpha^2_n}} = c^2_{n}.
\end{eqnarray}

Since we supposed $\|p\|_{L^2_w}=\sum_{i=0}^{n} c^2_i=1$, if $n=k$ then 
\begin{eqnarray}\label{MZGS}
\eta&=&\max_{p \in {\mathbb{P}}_n, \|p\|_{L^2_w}=1} |I(p^2)-S_m(p^2)| = \max_{p \in {\mathbb{P}}_n, \|p\|_{L^2_w}=1} {\color{black}c^2_n} = 1.
\end{eqnarray}

Hence for $m$-points Gaussian rules, when $n=k$ the quantity $\eta$ is $1$. Since $\eta$ is non decreasing in $n$, we have that $\eta \geq 1$ for $n \geq k$.

Thus the following theorem holds
\begin{theorem}
Let $S_k$ be a k-points Gaussian rule w.r.t. the weight function $w:(a,b) \rightarrow \mathbb R^+$, with $(a,b)$ not necessarily bounded. Then for $n \geq k$ we have
\begin{equation}\label{MZmax2b}
\eta=\max_{p \in {\mathbb{P}}_n, \|p\|_{L^2}=1} |I(p^2)-S(p^2)| \geq 1.
\end{equation}
\end{theorem}

\begin{Remark}
Observe that from (\ref{MZGS}) we also get that the worst case polynomials are $\hat{p}_n=\pm \psi_n$.
\end{Remark}

\begin{Remark}
{\color{black}It is common, in view of its simplicity, to use tensor-product rules $S_k=S_k^{(1)}\otimes S_k^{(2)}$ based on univariate Gaussian formulas with $m=2k-1$, to compute double integrals with respect to a product measure $d\mu=d\mu_1\,d\mu_2=w_1(x)w_2(y)\,dxdy$ on a box $(a,b)\times (c,d)$ (see {\cite[p.25]{ST71}}). In such cases, taking $n=k$ and the polynomial $p(x,y)=\psi_n(x)/\sqrt{\mu_2((c,d))}$, it is easily checked that $|I(p^2)-S_k(p^2)|=|1-S_k^{(1)}(\psi_n^2)|=1$. Hence, in view of the fact that $\eta$ is nondecreasing in $n$, if $n \geq k$ then $\eta \geq 1$. The same clearly holds true for product measures and tensor-product Gaussian rules on boxes in any dimension.
}
\end{Remark}

\section{Numerical examples}
In this section we intend to evaluate the Marcinkiewicz-Zygmund constants of some popular  rules over relevant sets such as the interval, the square, the disk, the triangle and the sphere.
{\color{black}In such cases we can refer to unit reference sets, since MZ inequalities of cubature rules for absolutely continuous measures, as well as conditioning of the corresponding Gramians, are invariant under affine transformations.} The open-source Matlab routines accompanying this paper are available at {\cite{S25}.

\subsection{The interval}
In the case of the interval $[-1,1]$, the Gauss-Legendre and the Clenshaw-Curtis rules are particularly appealing for numerical integration w.r.t. {\color{black}the Lebesgue measure}, cf. {\cite{TFT}}.
We have implemented a Matlab routine that for a fixed algebraic degree of {\color{black}exactness} $m$ and polynomial degree $n$: 
\begin{itemize}
{\color{black}
\item evaluates by three terms recursion the first $d_n=n+1$ orthonormal Legendre polynomials at the quadrature nodes;
\item computes for each rule with $ADE=m$ the Gramian $G$ and then defines the matrix $E=Id-G$;
\item determines the Marcinkiewicz-Zygmund constant $\eta=\|E\|_2$.
}
\end{itemize}

\begin{figure}[ht!]
\centering
\includegraphics[scale=0.42,valign=c]{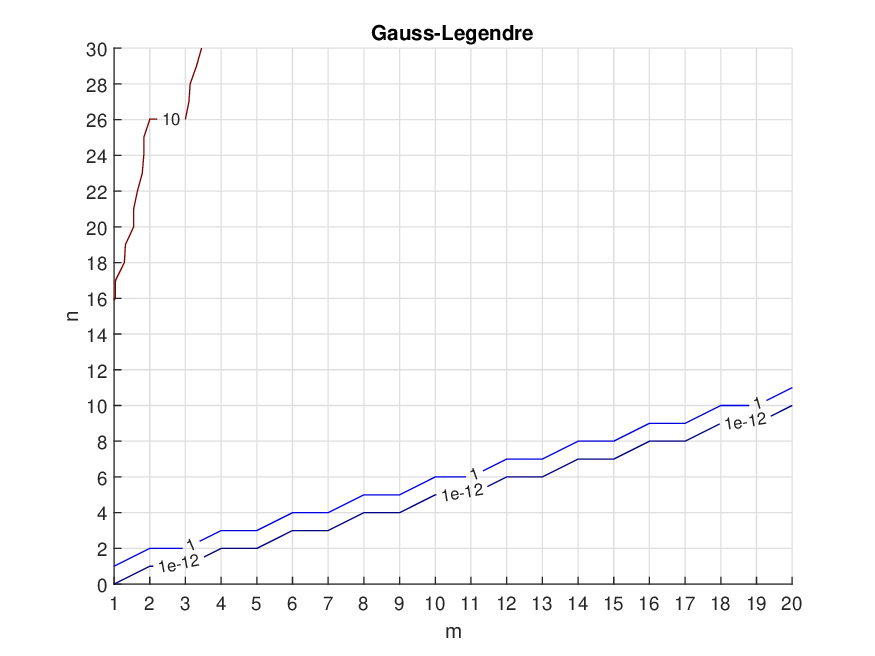} {\hspace{0.3cm}}
\includegraphics[scale=0.42,valign=c]{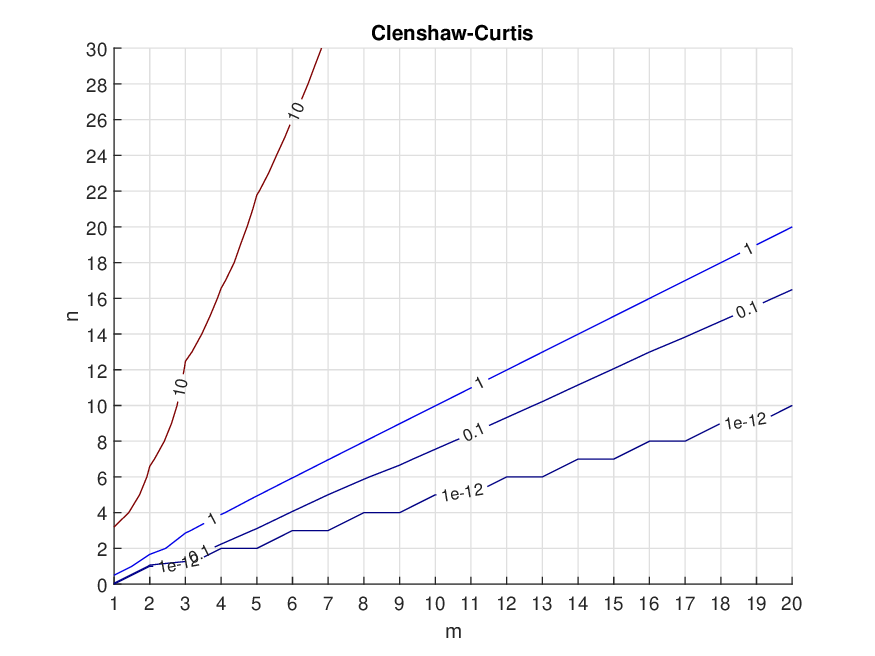} 
\caption{{\color{black}Contour lines of Marcinkiewicz-Zygmund constants $\eta$ in (\ref{etabound}) for Gauss-Legendre (left) and Clenshaw-Curtis (right) rules with $ADE=m$, where $m=1,2,\ldots,20$ and $n=0,\ldots,30$.}}
\label{MZ_1D_contour}
\end{figure}

In Figure \ref{MZ_1D_contour} we plot the contour lines of the Marcinkiewicz-Zygmund constants $\eta$ of Gauss-Legendre and Clenshaw-Curtis rules with $ADE=m$ for $m=1,\ldots,20$ and $n=0,\ldots,30$. 
In the Gaussian case it is clear that there is a jump when $n$ passes from $k=\lfloor\frac{m}{2}\rfloor$ to $k+1$ since the MZ constant varies abruptly from values close to machine precision to $1$, as foreseen by the theory.
{\color{black}The variation is slower in the case of the Clenshaw-Curtis rule, where the contour lines $\eta\approx 0$ and $\eta=1$ diverge  significantly, and for a given $m$ the MZ constant $\eta$ turns out to be strictly 
less than 1 for $\lfloor\frac{m}{2}\rfloor<n\leq m-1$. Such a behavior, which to our knowledge is not known, could deserve further study. It could be related to the fact, manifest also in Figure 1, that Clenshaw-Curtis rule is able to reasonably approximate  integrals of polynomials for some degree higher than what is expected, a phenomenon already known and studied in \cite{TFT}.

About the conditioning ${\mbox{cond}}_2(G)$ of the Gramian matrix based on the two rules above (computed by the Matlab function {\tt cond}), the behaviour is also different and is illustrated in Figure \ref{MZ_1D_cond}. For the 
\begin{itemize}
\item Gauss-Legendre rule  ${\mbox{cond}}_2(G) \approx 1$ when $n \leq \lfloor m/2 \rfloor$ otherwise is larger than the inverse of machine precision,
\item Clenshaw-Curtis ${\mbox{cond}}_2(G)$ is smaller than $10$ for $\lfloor m/2 \rfloor <n\leq m\;,$ whereas it increases 
abruptly up to near-singularity of $G$ for $n>m$.
\end{itemize}
}

\begin{figure}[ht!]
\centering
\includegraphics[scale=0.42,valign=c]{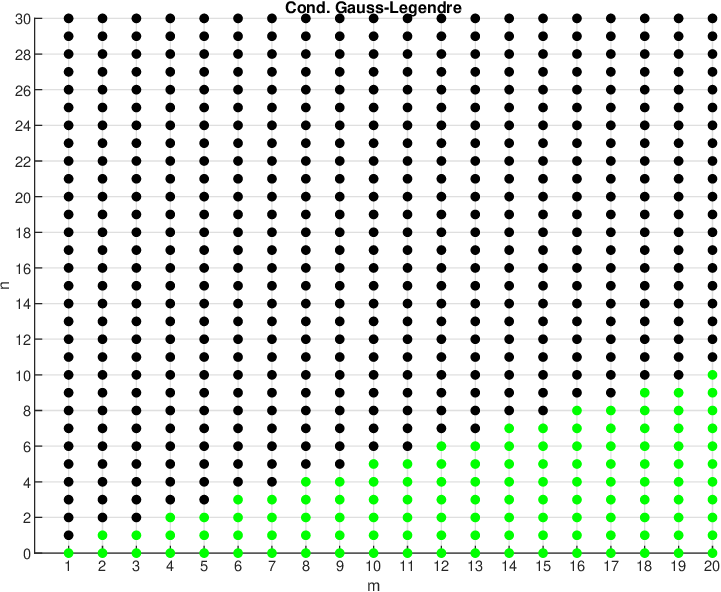} 
\quad
\includegraphics[scale=0.42,valign=c]{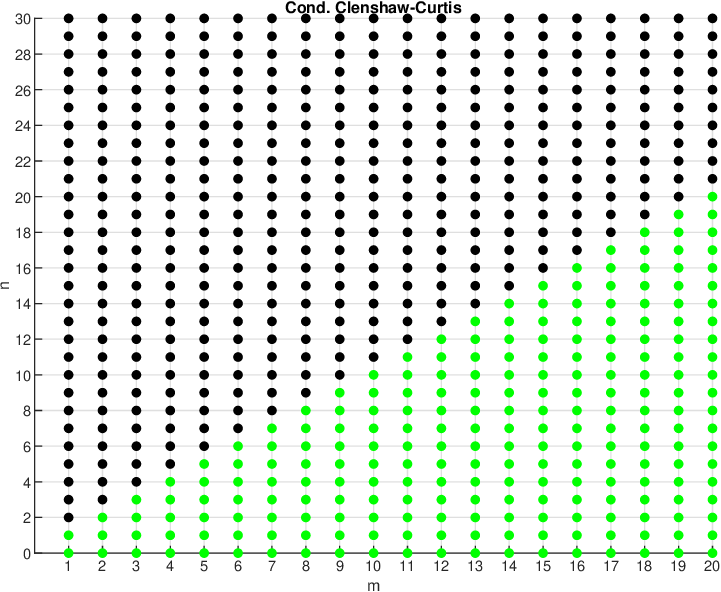} 
\caption{{\color{black}Conditioning in norm 2 of the Gramian based on Gauss-Legendre and Clenshaw-Curtis rules with $ADE=m$, where $m=1,2,\ldots,20$ and $n=0,\ldots,30$. 
Green dot: $\mbox{cond}_2(G) \in [1,10)$. 
Black dot: $\mbox{cond}_2(G)\in [10^7,+\infty)$.}}
\label{MZ_1D_cond}
\end{figure}

\subsection{The square}

Here we compute the MZ constants of some rules on the unit-square $[-1,1]^2$ w.r.t. {\color{black}the Lebesgue measure and the product Chebyshev measure} (see \cite{DR77,ST71}). A first example is the tensor product rule based on Gauss-Legendre rules. 
Alternative rules for the Lebesgue measure are that based on Padua points \cite{CDMSV} and the near minimal rules proposed in \cite{FS12}. As  last example we take into account the {\it{Morrow-Patterson-Xu}} cubature rule, that is (near) minimal for the product Chebyshev measure (see \cite{MP,X96}). 
{\color{black}
We first compute at the cubature nodes the total-degree orthonormal product basis of dimension $d_n=(n+1)(n+2)/2$ for the given measure, then the error matrix $E=Id-G$, where $G$ is the corresponding Gramian, and finally $\eta=\|E\|_2$. 

The numerical results are collected in Figures (\ref{figure:MZ_square})-(\ref{figure:MZ_square_cond}). 
We notice that in the case of the tensor-product, near minimal  and Morrow-Patterson-Xu rules, for $m=1,\ldots,20$ and $n > \lfloor m/2 \rfloor$ it is always $\eta \geq 1$, and  
$\mbox{cond}_2(G)$ increases abruptly up to near-singularity of $G$. Differently, 
with Padua points based cubature, $\eta < 1$ up to $n=m-1$, and $\mbox{cond}_2(G)$ is smaller than $10$ for $n<m\;,$ whereas it increases 
abruptly up to near-singularity of $G$ for $n>m$.
}

\begin{figure}[ht!]
\centering
\includegraphics[scale=0.42,valign=c]{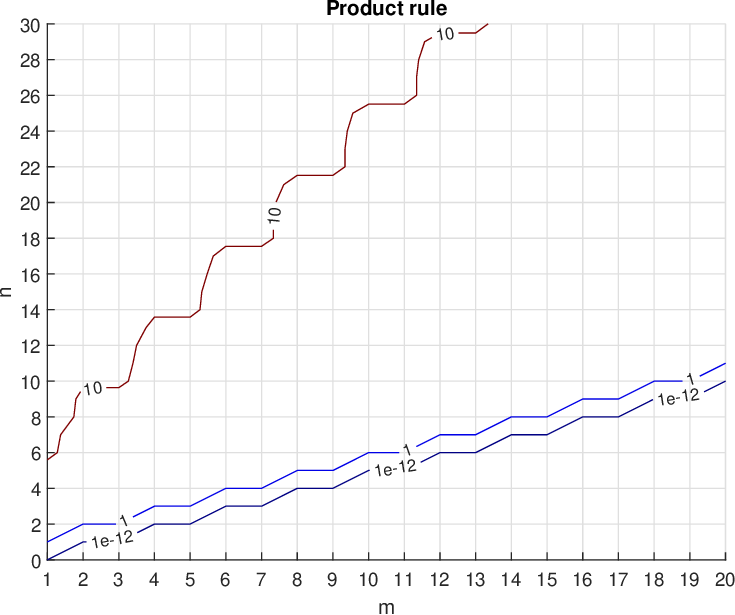}
{\hspace{0.3cm}}
\includegraphics[scale=0.42,valign=c]{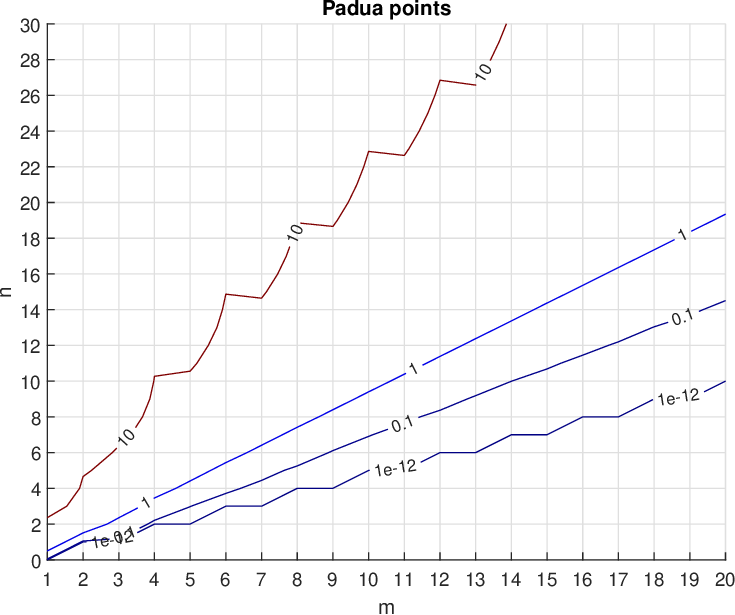} \\
\includegraphics[scale=0.42,valign=c]{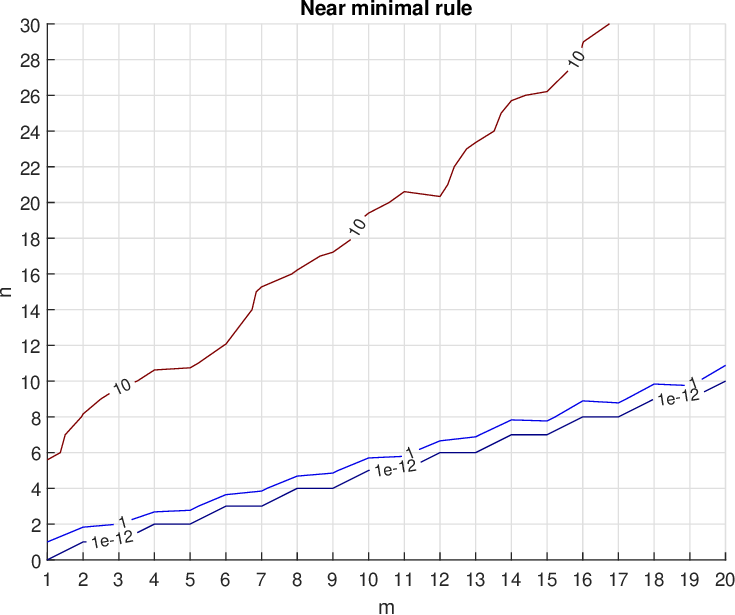}
{\hspace{0.3cm}}
\includegraphics[scale=0.42,valign=c]{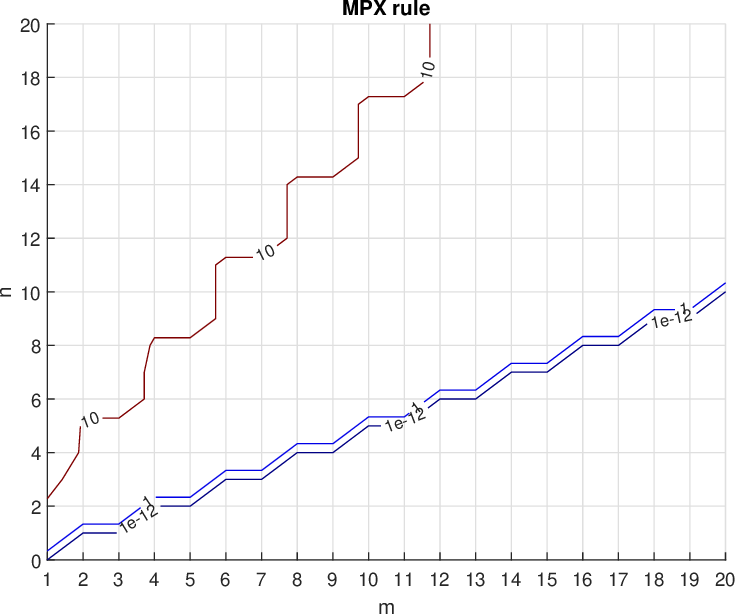}
\caption{{\color{black}Contour lines of Marcinkiewicz-Zygmund constants $\eta$ in (\ref{etabound}) some cubature rules on the unit-square. Tensor product Gauss-Legendre rule (top-left), Padua points based rule (top-right), near-minimal (bottom-left) and Morrow-Patterson-Xu (bottom-right) rules.}}
\label{figure:MZ_square}
\end{figure}

\begin{figure}[ht!]
\centering
\includegraphics[scale=0.42,valign=c]{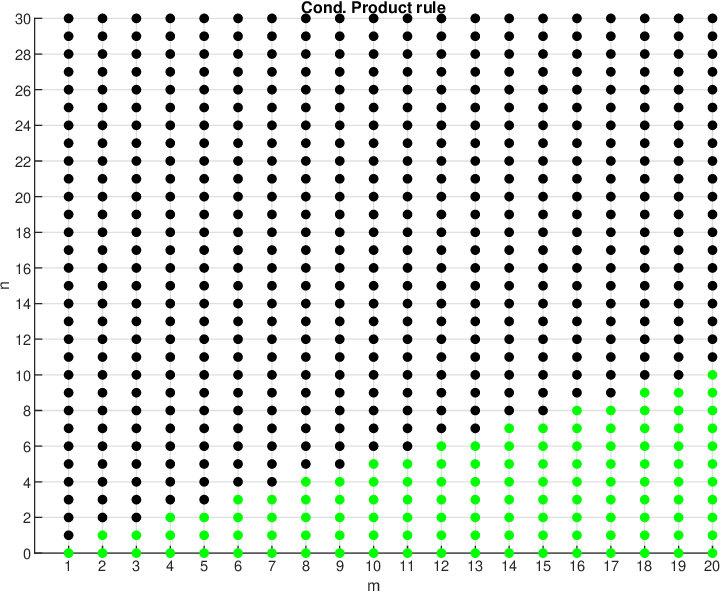} 
\quad
\includegraphics[scale=0.42,valign=c]{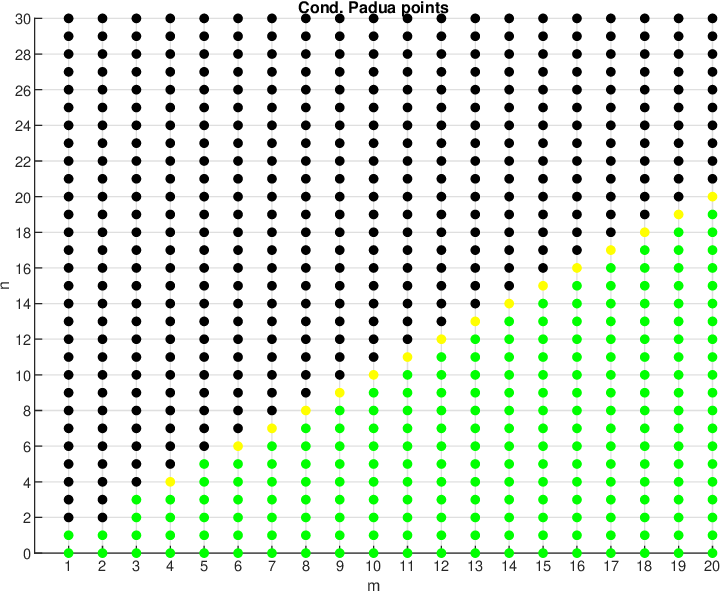} \\
\includegraphics[scale=0.42,valign=c]{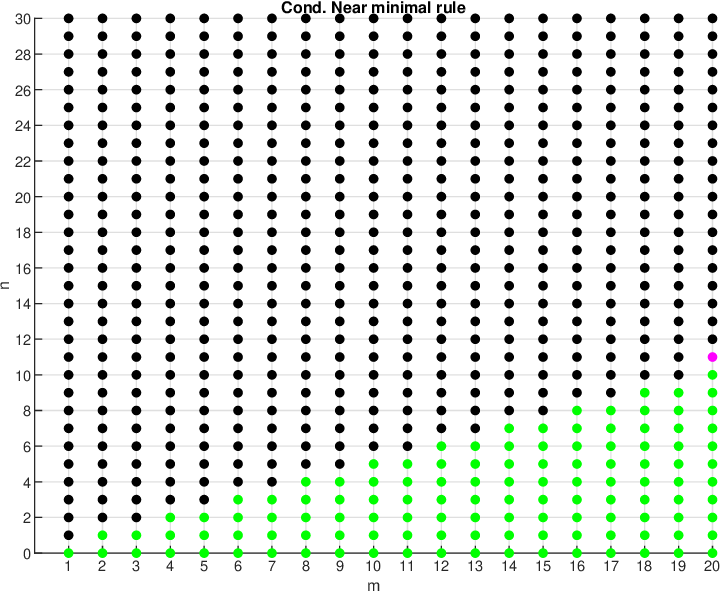} 
\quad
\includegraphics[scale=0.42,valign=c]{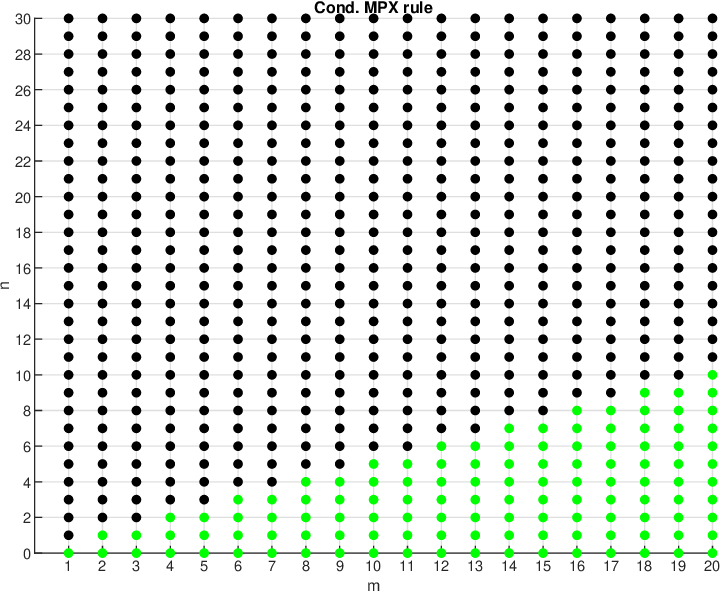}
\caption{{\color{black}Conditioning in norm 2 of the Gramian based on tensor-product Gauss-Legendre rule, Padua points based rule, near minimal rule and and Morrow-Patterson-Xu rule with $ADE=m$, where $m=1,2,\ldots,20$ and $n=0,\ldots,30$. Green dot: $\mbox{cond}_2(G) \in [1,10)$. Yellow dot: $\mbox{cond}_2(G) \in [10,10^ 2)$. Magenta dot: $\mbox{cond}_2(G)\in [10^4,10^7)$.
Black dot: $\mbox{cond}_2(G)\in [10^7,+\infty)$.}
}
\label{figure:MZ_square_cond}
\end{figure}
 
{\color{black}
\subsection{The disk and the triangle}
In the case of the unit-disk $B(0,1)=\{(x,y): x^2+y^2 \leq 1\}$ we have used the Logan-Shepp orthonormal basis w.r.t. the Lebesgue measure \cite{LS75}. tensor-product rules in polar coordinates are a popular choice, based on Gauss-Legendre and trapezoidal rules, with $ADE=m$ (see {\cite[p.32]{ST71}}). The advantage of these rules is that they are easily implemented even for large $m$, though their cardinality is far from being minimal.
Alternatively, for mild degrees $m$ there are rules that have $ADE=m$ but cardinality particularly low (see, e.g. \cite{S25sets}). For small $m$ they have been determined theoretically so that the number of nodes is minimal while for higher values they have been  obtained numerically. 

In the case of the unit-simplex 
$T=\{(x,y): x,y\geq 0,x+y\leq 1\}$ we have used the Dubiner orthonormal basis w.r.t. the Lebesgue measure \cite{Dub91}. Tensor-product rules are again a popular choice, being easily implemented at any degree of precision (see, e.g. {\cite[p.29]{ST71}} and (\cite{LC93}). They are generally named {\it{Stroud Conical Rules}}, have positive weights, internal nodes, $ADE=m$ on the reference triangle, with cardinality $(\lceil\frac{m+1}{2} \rceil)^2 \approx \frac{m^2}{4}$.
Again, for mild degrees, there are rules that have $ADE=m$, but a lower cardinality. Some of them are known theoretically, while some others via optimization algorithms (see {\cite{LJ75}} and {\cite{D85}}). We refer to {\cite{S25sets}} for recent updates of these formulas. 

The numerical results are collected in Figures (\ref{figure:MZ_disk})-(\ref{MZ_simplex_cond}). Only for the tensor-product rule $\eta<1$ for $n>\lfloor\frac{m}{2}\rfloor$, since the contour lines $\eta\approx 0$ and $\eta=1$ slightly diverge. On the other hand,    
in all cases $\mbox{cond}_2(G)$ increases abruptly up to near-singularity of $G$ for $n > \lfloor m/2 \rfloor$.
}

\begin{figure}[ht!]
\centering
\includegraphics[scale=0.42,valign=c]{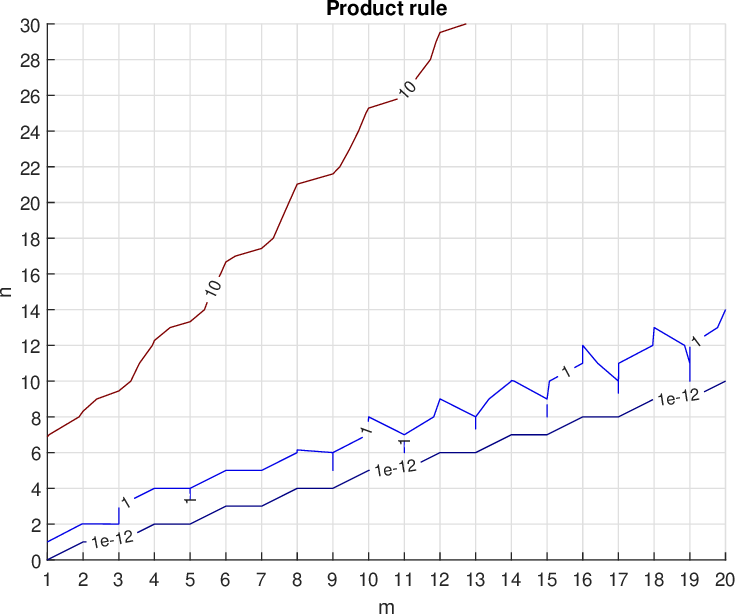} {\hspace{0.3cm}}
\includegraphics[scale=0.42,valign=c]{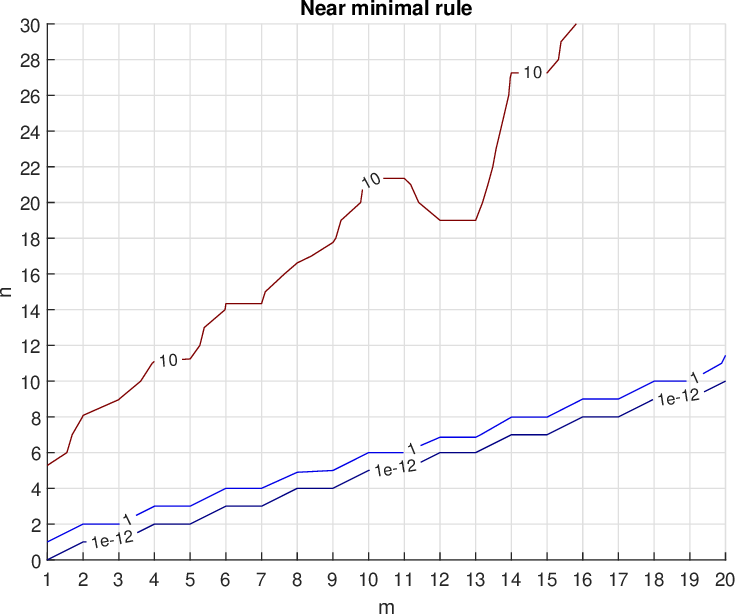}
\caption{{\color{black}Contour lines of Marcinkiewicz-Zygmund constants $\eta$ in (\ref{etabound}) for a couple of cubature rules on the unit-disk.}}
\label{figure:MZ_disk}
\end{figure}

\begin{figure}[ht!]
\centering
\includegraphics[scale=0.42,valign=c]{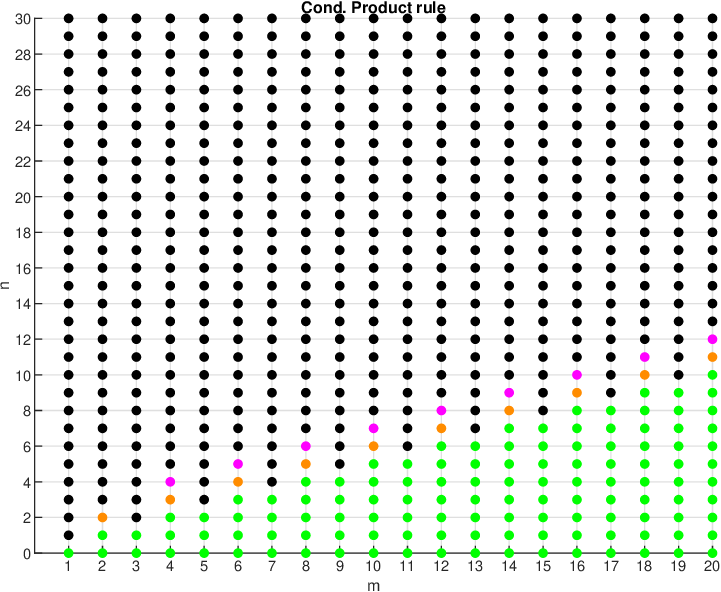} 
\quad 
\includegraphics[scale=0.42,valign=c]{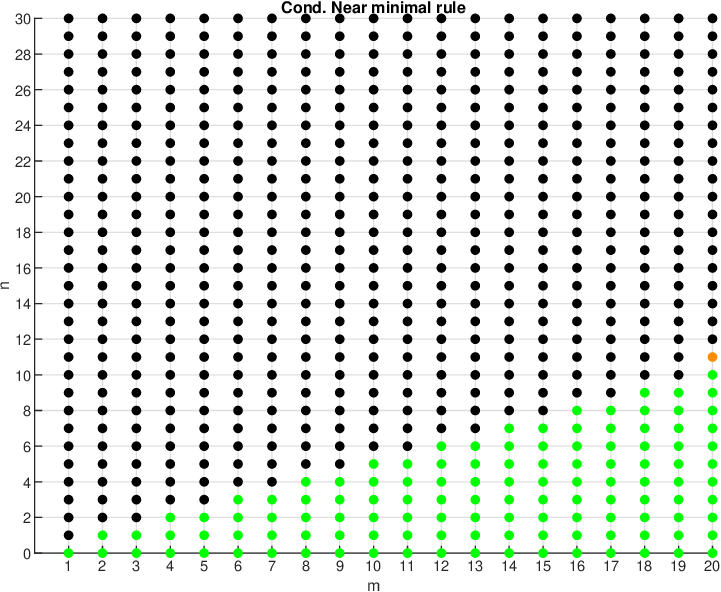} 
\caption{{\color{black}Conditioning in norm 2 of the Gramian for a couple of cubature rules on the unit-disk. 
Green dot: $\mbox{cond}_2(G) \in [1,10)$. Orange dot: $\mbox{cond}_2(G) \in [10^2,10^4)$. Magenta dot: $\mbox{cond}_2(G)\in [10^4,10^7)$.
Black dot: $\mbox{cond}_2(G)\in [10^7,+\infty)$.}}
\label{MZ_disk_cond}
\end{figure}

\begin{figure}[ht!]
\centering
\includegraphics[scale=0.42,valign=c]{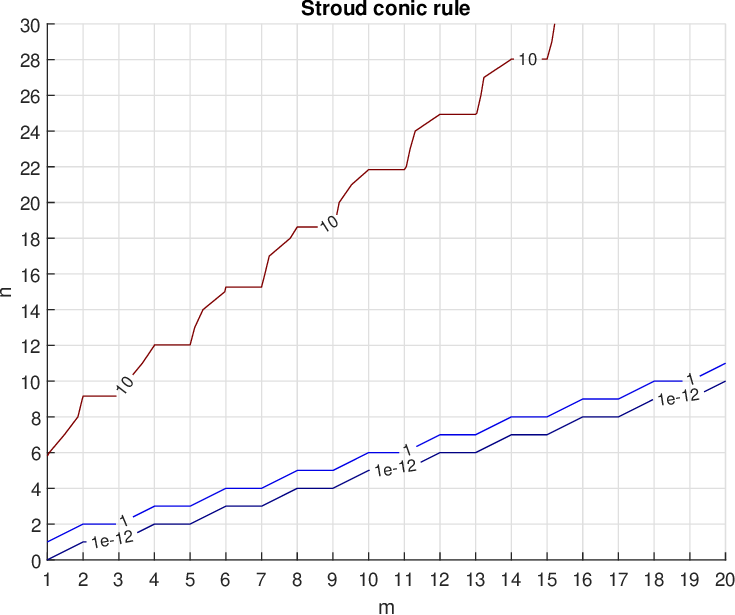}
\quad
\includegraphics[scale=0.42,valign=c]{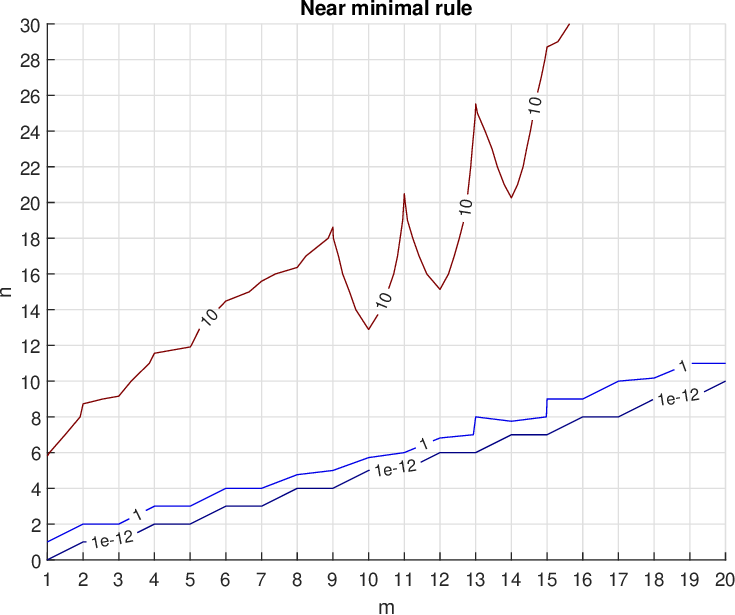} 
\caption{{\color{black}Contour lines of Marcinkiewicz-Zygmund constants $\eta$ in (\ref{etabound}) for a couple of cubature rules on the unit-simplex.}}
\label{figure:MZ_simplex}
\end{figure}

\begin{figure}[ht!]
\centering
\includegraphics[scale=0.42,valign=c]{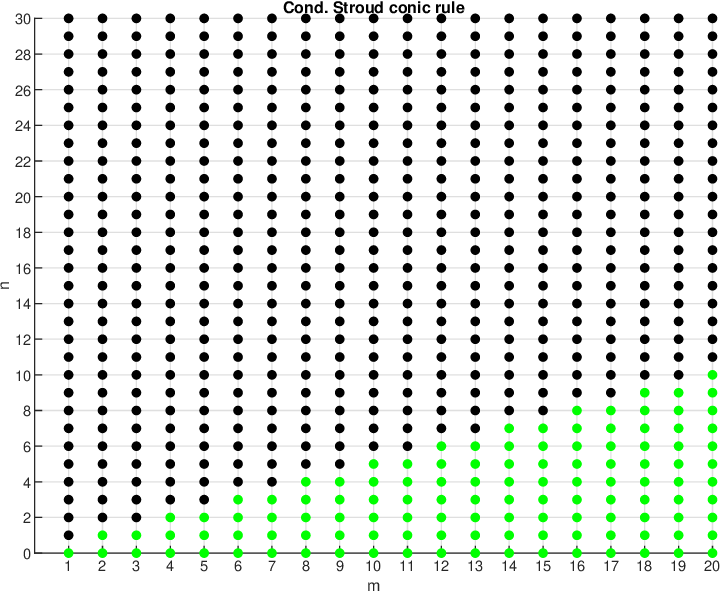} 
\quad
\includegraphics[scale=0.42,valign=c]{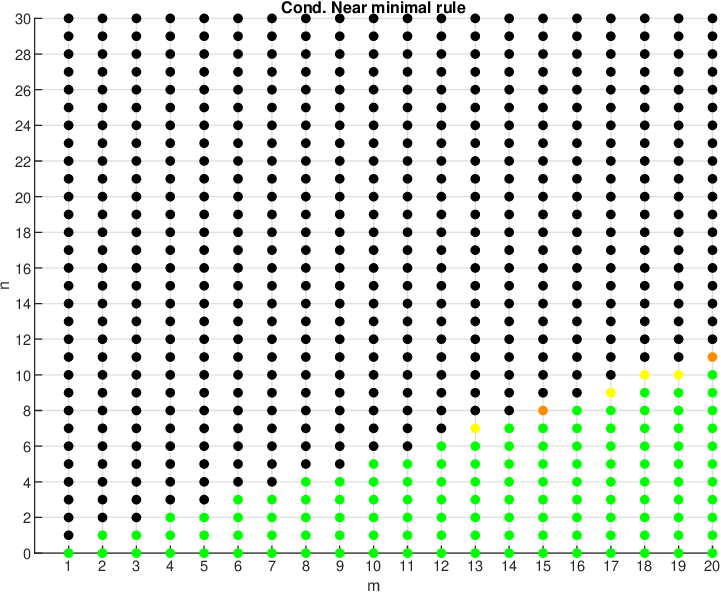} 
\caption{{\color{black}Conditioning in norm 2 of the Gramian 
for a couple of cubature rules on the unit-simplex.
Green dot: $\mbox{cond}_2(G) \in [1,10)$. Yellow dot: $\mbox{cond}_2(G) \in [10,10^2)$. Orange dot: $\mbox{cond}_2(G) \in [10^2,10^4)$. 
Black dot: $\mbox{cond}_2(G)\in [10^7,+\infty)$.}}
\label{MZ_simplex_cond}
\end{figure}

\subsection{The sphere}

On the sphere, in addition to the basic {\it{latitude-longitude}} rule obtained by reformulating the problem in spherical coordinates and applying suitable Gaussian and trapezoidal rules {\cite[p.40]{ST71}}, there are many works concerning spherical designs (see, e.g., the pioneering paper \cite{DGS77} as well as \cite{HSW15} for an introduction on the topic). Here we focus our attention to the recent rules proposed in {\cite{RWH}}, that are a specific set of {\it{spherical t-designs on ${\mathbb{S}}^2$}}, for $m=1,2,\ldots,180$ and  {\it{symmetric (antipodal) spherical t-designs on ${\mathbb{S}}^2$}} for $m=1,3,5,\ldots,325$. These Matlab datasets are available as open-source codes.

In our numerical tests, we adopted as orthonormal basis that of spherical harmonics. The numerical tests are illustrated in Figure {\ref{figure:MZ_sphere}}. 
While for the {\it{latitude-longitude}} rule the Marcinkiewicz-Zygmund constant $\eta$ is smaller than $1$ only when $n \leq \lfloor m/2 \rfloor$,  the spherical designs may have $\eta < 1$ even for some $n > \lfloor m/2 \rfloor$,  {\color{black}since the contour lines $\eta\approx 0$ and $\eta=1$ slightly diverge. On the other hand, also the Gramians of the spherical designs remain well-conditioned up to $n$ slightly larger than $ \lfloor m/2 \rfloor$ and then the conditioning increases abruptly.}

\begin{figure}[ht!]
\centering
\includegraphics[scale=0.42,valign=c]{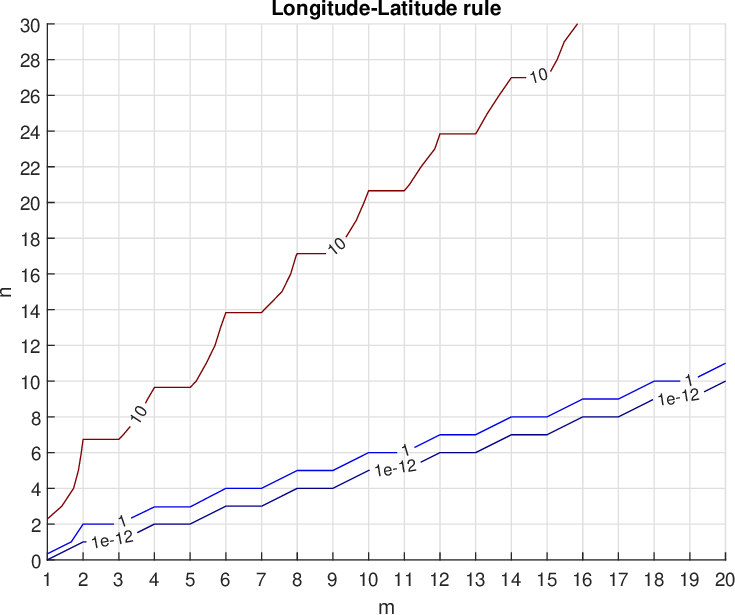}{\hspace{0.1cm}}
\includegraphics[scale=0.42,valign=c]{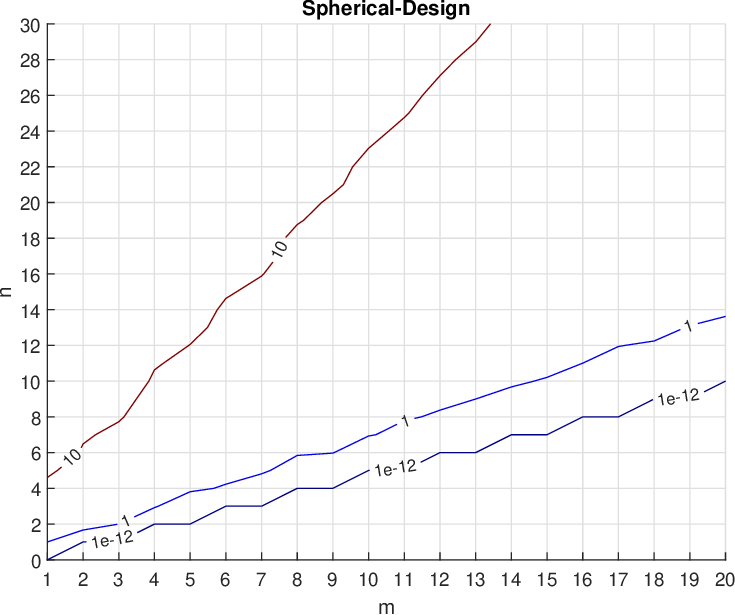}{\hspace{0.1cm}}
\includegraphics[scale=0.42,valign=c]{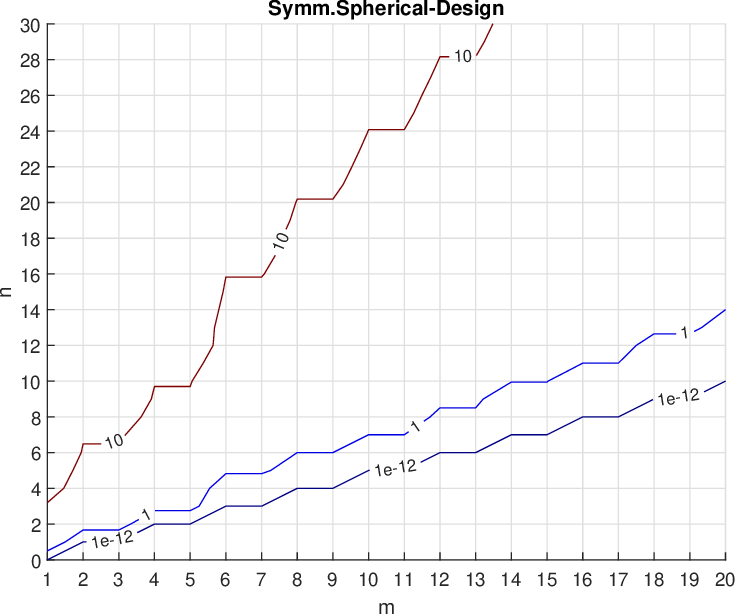}
\caption{{\color{black}Contour lines of Marcinkiewicz-Zygmund constants $\eta$ in (\ref{etabound}) for some cubature rules on the unit-sphere ${\mathbb{S}}^2$. From left to right, latitude-longitude rules, spherical designs and symmetrical spherical designs suggested in \cite{RWH}.}}
\label{figure:MZ_sphere}
\end{figure}

\begin{figure}[ht!]
\centering
\includegraphics[scale=0.42,valign=c]{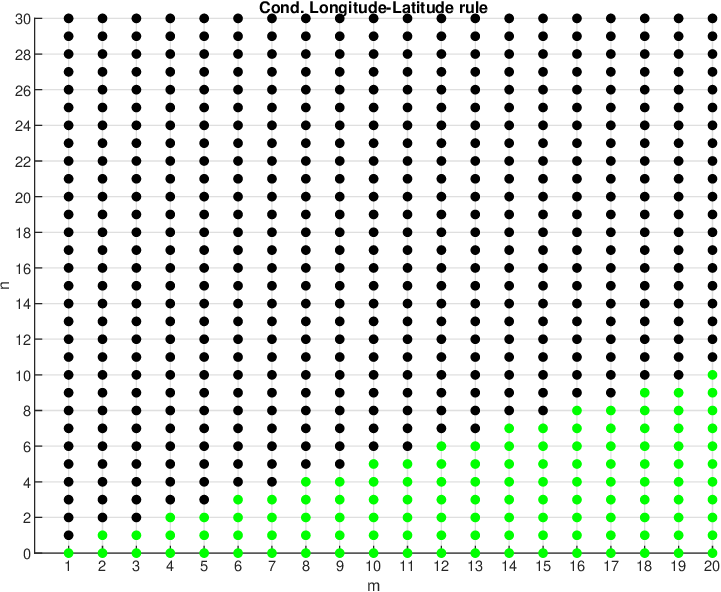}
{\hspace{0.1cm}}
\includegraphics[scale=0.42,valign=c]{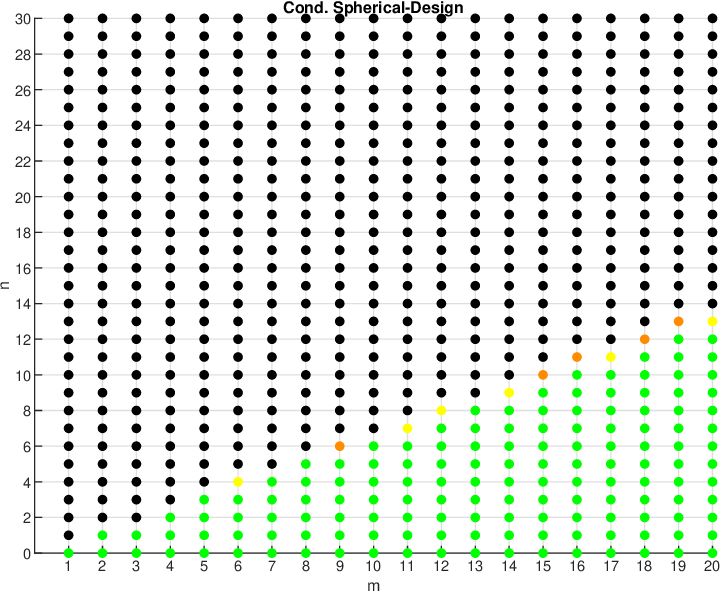}
{\hspace{0.1cm}}
\includegraphics[scale=0.42,valign=c]{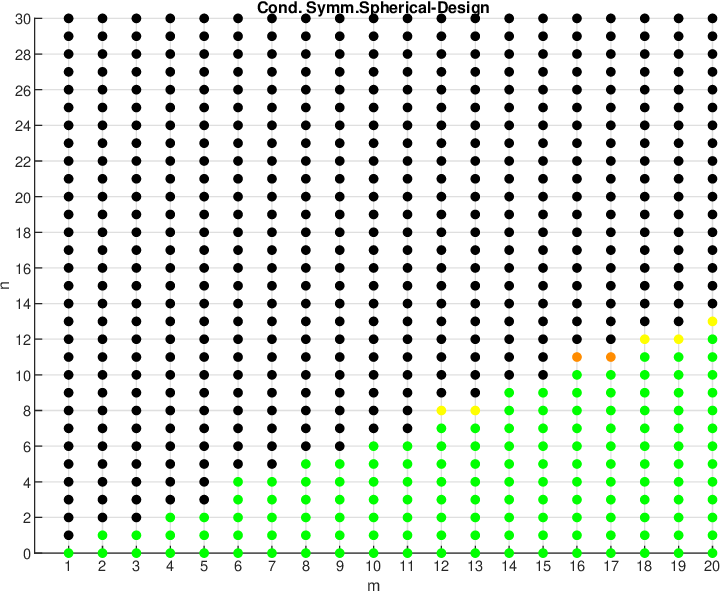} 
\caption{{\color{black}Conditioning in norm 2 of the Gramian for some cubature rules on the unit-sphere $\mathbb{S}^2$.
Green dot: $\mbox{cond}_2(G) \in [1,10)$. 
Yellow dot: $\mbox{cond}_2(G) \in [10,10^2)$. Orange dot: $\mbox{cond}_2(G) \in [10^2,10^4)$. Magenta dot: $\mbox{cond}_2(G)\in [10^4,10^7)$.
Black dot: $\mbox{cond}_2(G)\in [10^7,+\infty)$.}}
\label{MZ_sphere_cond}
\end{figure}

\subsection{QMC rules}

We conclude our analysis on cubature formulas paying attention to QMC rules. We have tested on the unit-square $[-1,1]^2$ and on the unit-cube $[-1,1]^3$ QMC rules based on Halton sets. {\color{black}Since the reference measure is the Lebesgue measure, we have chosen as underlying orthonormal basis of $\mathbb{P}_n$ the total-degree product Legendre basis.} 
We have considered $M=2^{m}$, $m=1,\ldots,20$ as cardinality of the rule, and $n=0,\ldots,20$.

The respective results on these domains about the values of $\eta$ and the conditioning in norm $2$ of the Gramian are illustrated, respectively, in Figures \ref{figure:MZ_square_QMC} and \ref{figure:MZ_cube_QMC}.
As suggested by the convergence of the Gramian to the identity matrix as $M\to \infty$, for a fixed $n$, $\eta$ decreases while increasing $M$. The higher is $M=2^m$ the larger are the values of $n$ for which $\eta<1$. 
We observe that there are couples $(m,n)$ for  which $\eta > 1$, but ${\mbox{cond}}(G)$ is not too large. This implies that we can still {\it{safely}} solve the Least-Squares problem in $\mathbb{P}_n$ via Gramian matrix. {\color{black}Moreover, the first estimate in (\ref{alt}) holds, whereas the bound (\ref{QMC-err}) does not apply.}

\begin{figure}[ht!]
\centering
\includegraphics[scale=0.42,valign=c]{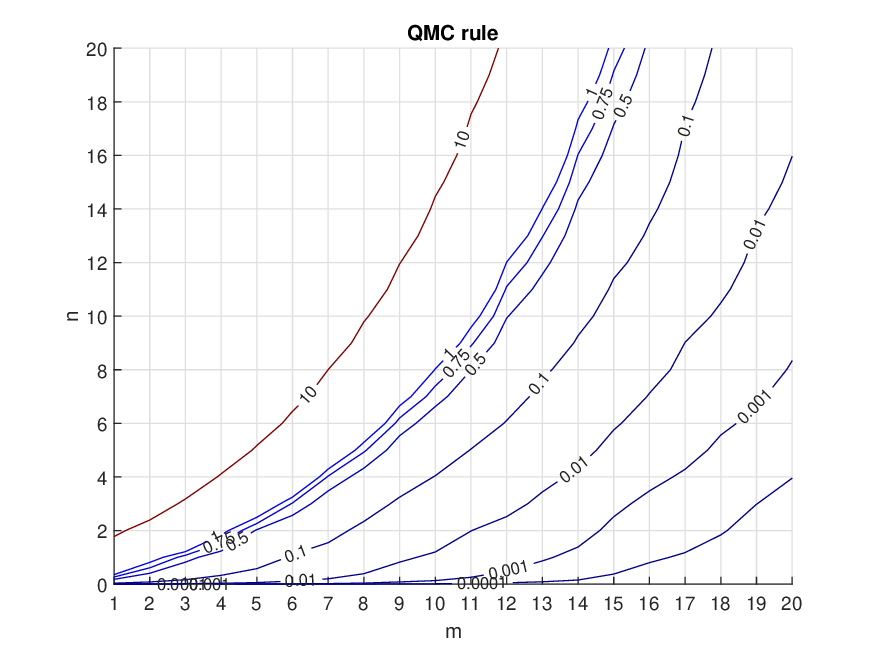}
\quad
\includegraphics[scale=0.42,valign=c]{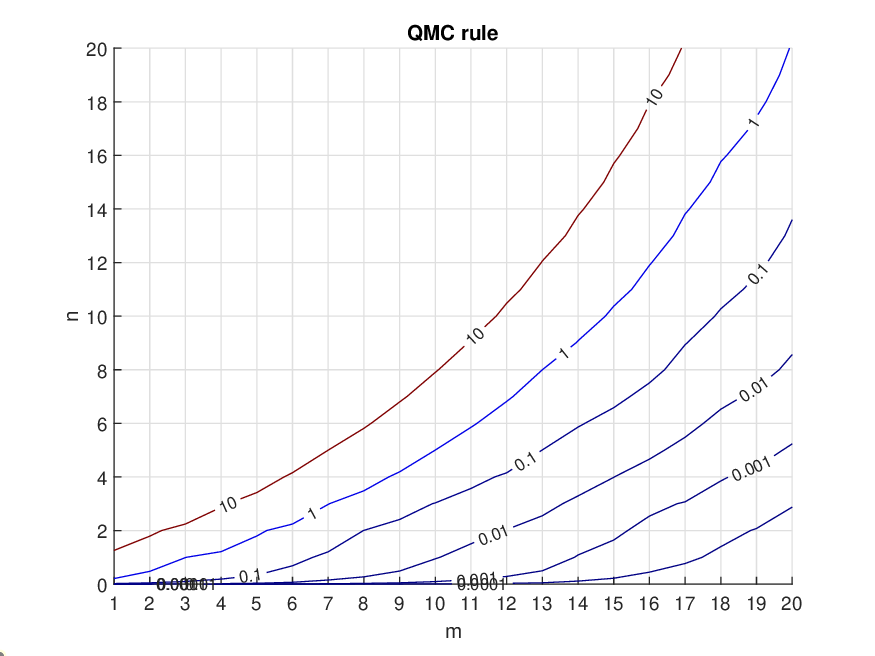}

\caption{On the left: Contour lines of Marcinkiewicz-Zygmund constants $\eta$ in (\ref{etabound}) for QMC rules on the unit-square $[-1,1]^2$. On the right: Contour lines of Marcinkiewicz-Zygmund contants of QMC rules in (\ref{etabound}) on the unit-cube $[-1,1]^3$. In this example we consider QMC rules at Halton points in which the cardinality is $M=2^m$ for $m=1,\ldots,20$, and $n=0,\ldots,20$.}
\label{figure:MZ_square_QMC}
\end{figure}

\begin{figure}[ht!]
\centering
\includegraphics[scale=0.42,valign=c]{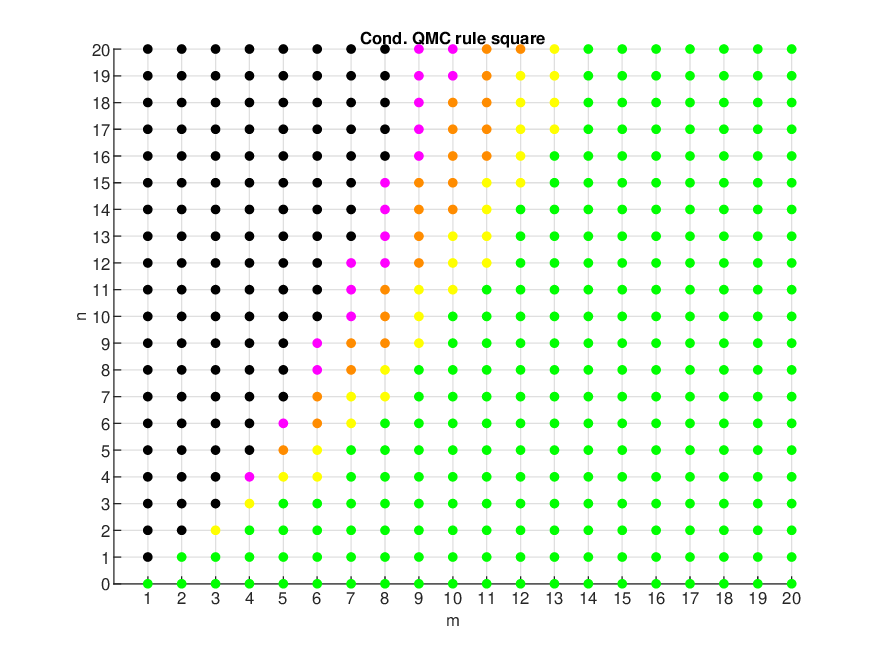}
\quad
\includegraphics[scale=0.42,valign=c]{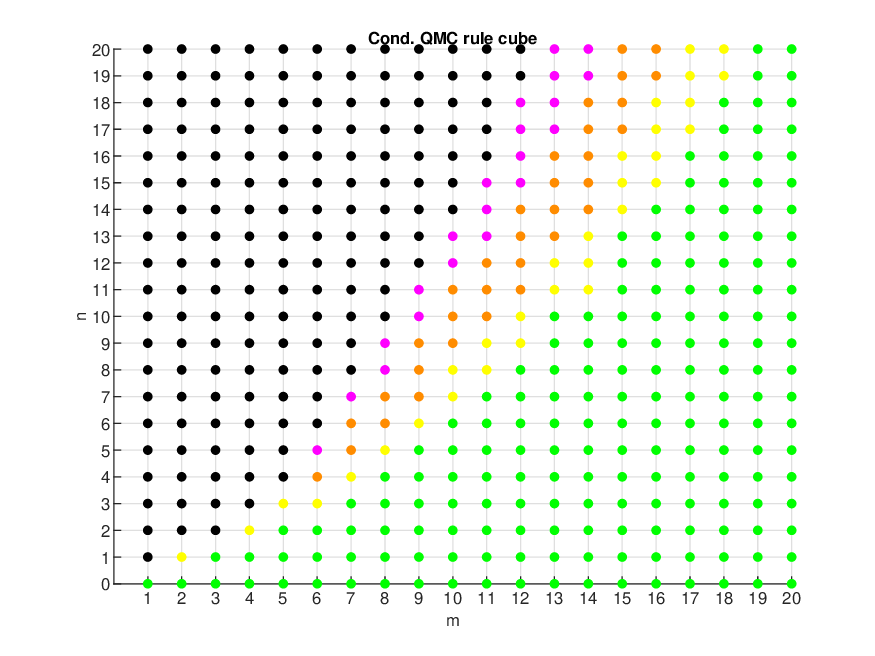}
\caption{On the left: Conditioning in norm 2 of the Gramian  based on QMC rules on the unit-square $[-1,1]^2$. On the right: Conditioning in norm 2 of the Gramian  based on QMC rules on the unit-cube $[-1,1]^3$. Green dot: $\mbox{cond}_2(G) \in [1,10)$. 
Yellow dot: $\mbox{cond}_2(G) \in [10,10^2)$. Orange dot: $\mbox{cond}_2(G) \in [10^2,10^4)$. Magenta dot: $\mbox{cond}_2(G)\in [10^4,10^7)$.
Black dot: $\mbox{cond}_2(G)\in [10^7,+\infty)$}
\label{figure:MZ_cube_QMC}
\end{figure}

{\color{black}
\subsection{Polynomial approximation by orthogonal bases}
Among the tested algebraic rules on different basic domains, Clenshaw-Curtis rule for the interval and Padua points based rule for the square clearly emerged as good candidates for exactness-relaxed (or unfettered) hyperinterpolation as proposed in \cite{ANWU22}. 
We guess that this could be due to special properties of the Chebyshev polynomials, on the line of the study in \cite{TFT} and might deserve further investigation, that however goes beyond the scope of the present paper. 
Moreover, following \cite{AKG25}, QMC rules with sufficiently large cardinality $M$ can also be used for unfettered hyperinterpolation in the square and the cube. The reason is that $A,B\to 1$ as $M\to \infty$ for fixed $n$, since the Gramian converges to the identity matrix, so $\eta\to 0$ and thus is bounded away from 1 for $M=M(n)$ sufficiently large. 

It is then worth comparing, also in view of estimates (\ref{QMC-err}) and (\ref{err-Gn}), the relative errors in the $L^2$-norm of classical Sloan's hyperinterpolation with both, unfettered hyperinterpolation and Least-Squares approximation, constructed by the three rules just quoted. Since the reference
measure is the Lebesgue measure, we have chosen as underlying orthonormal basis of $\mathbb{P}_n$ the total-degree product Legendre basis. 
In particular we test Clenshaw-Curtis rule for the interval, Padua points based rule for the square and QMC rule at Halton points for the cube. 

In the first two cases we fix an $ADE$, say $m=15$. Classical hyperinterpolation ($\eta=0$) is implemented by a (product) Gauss-Legendre rule with $ADE$ $2m=30$. Given a test function $f$ and a polynomial hyperinterpolant $p_n=\mathcal{H}_nf$ as in (\ref{hyper}) or a Least-Squares projection $p_n=\mathcal{G}_nf$ as in (\ref{orth}), for $n=1,2,\dots,15$, we approximate the relative errors in the $L^2$-norm by a (product) Gauss-Legendre rule of $ADE=50$ 
with nodes say $X=\{\mathbf{z}_k\}$ and weights say $u=\{u_k\}$, as 
\begin{equation}\label{L2re}
\frac{\|p_n-f\|_{L^2}}{\|f\|_{L^2}}\approx \frac{\|p_n-f\|_{\ell^2_{u}(X)}}{\|f\|_{\ell^2_{u}(X)}}=\frac{\left(\sum_{k=1}^{\nu_{50}}{u_k\,(p_n(\mathbf{z}_k)-f(\mathbf{z}_k))^2}\right)^{1/2}}{\left(\sum_{k=1}^{\nu_{50}}{u_k\,f^2(\mathbf{z}_k)}\right)^{1/2}}\;,
\end{equation}
where $\nu_{50}=26$ for the interval and 
$\nu_{50}=26^2$ for the square. 

As test functions with different regularity we adopted for the unit-interval $[-1,1]$ 
$$
f_1(x)=\exp(-x^2), \, f_2(x)=(0.5+x)^{15}, \, f_3(x)=\sin(\pi x), \, f_4(x)=|x-0.5|^{3}, \,f_5(x)=|x-0.5|^{7}.
$$
and for the unit-square $[-1,1]^2$ 
$$
f_1(x,y)=\exp(-(x^2+y^2)), \, f_2(x,y)=(0.5+x+0.1 y)^{15}, \, f_3(x,y)=\sin(\pi x + \pi y),
$$
$$
f_4(x,y)=d((x,y),(0.5,0.5))^{3}=((x-0.5)^2+(y-0.5)^2)^{3/2}, 
$$
$$
\,f_5(x,y)=d((x,y),(0.5,0.5))^{7}=((x-0.5)^2+(y-0.5)^2)^{7/2}.
$$.

The numerical results show that though unfettered hyperinterpolation offers decent errors for $f_4$ and $f_5$ at low $n$, it tends to deteriorate increasing $n$. This is particularly evident for the most regular functions, where the minimum error of unfettered hyperinterpolation stays around $10^{-5}$. Notice that it does not well reconstruct the polynomial $f_2$, as expected because it is not a projection on $\mathbb{P}_n$. 

We observe that the errors obtained by Padua points based LS approximation are close to those of classical hyperinterpolation, but using a rule with $ADE=15$ instead of $ADE=30$. This can be seen as another way of relaxing the strict quadrature exactness of classical hyperinterpolation. While for the interval the sampling cardinality is 16 for both methods, in the case of the square LS samples at only $d_{15}=16\times 17/2=136$ Padua points, whereas classical hyperinterpolation via product Gauss-Legendre rule samples at $(m+1)^2=(15+1)^2=256$ points. This means that Padua points based LS behave as we used a {\em minimal} cubature rule with $ADE=30$ for the Lebesgue measure 
on the square (which is however not known, cf. the recent monograph \cite{X25}). On the other hand, the advantage would be maintained even if we used classical hyperinterpolation based on the best known {\em near-minimal} rule with $ADE=30$, obtained in {\cite{FS12}}, that has $167$ points.
}

\begin{figure}[ht!]
\centering
\includegraphics[scale=0.42,valign=c]{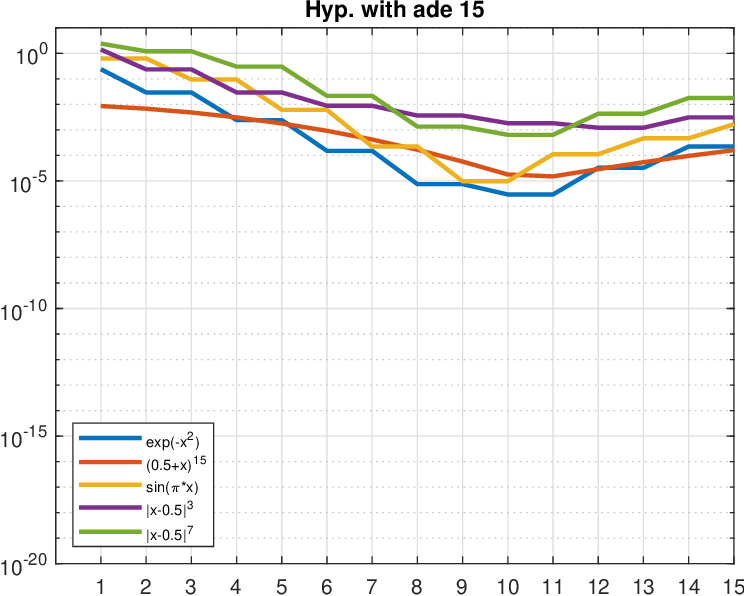}
\includegraphics[scale=0.42,valign=c]{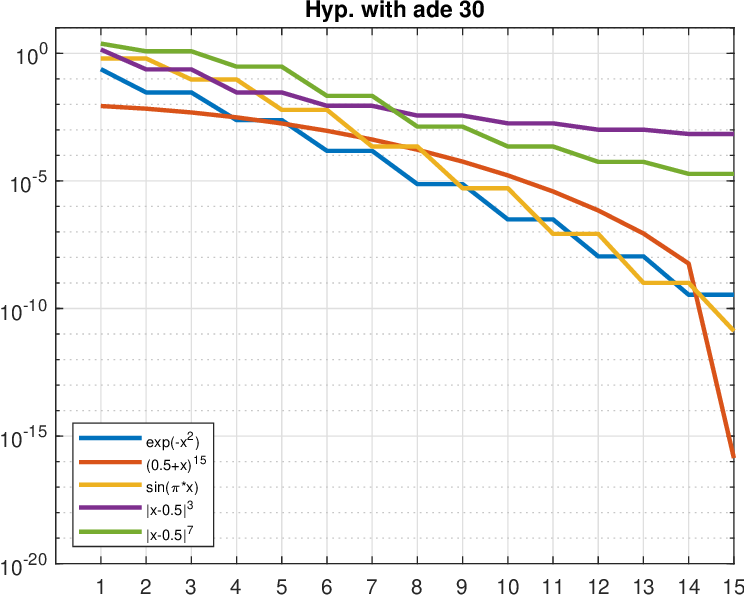}
\includegraphics[scale=0.42,valign=c]{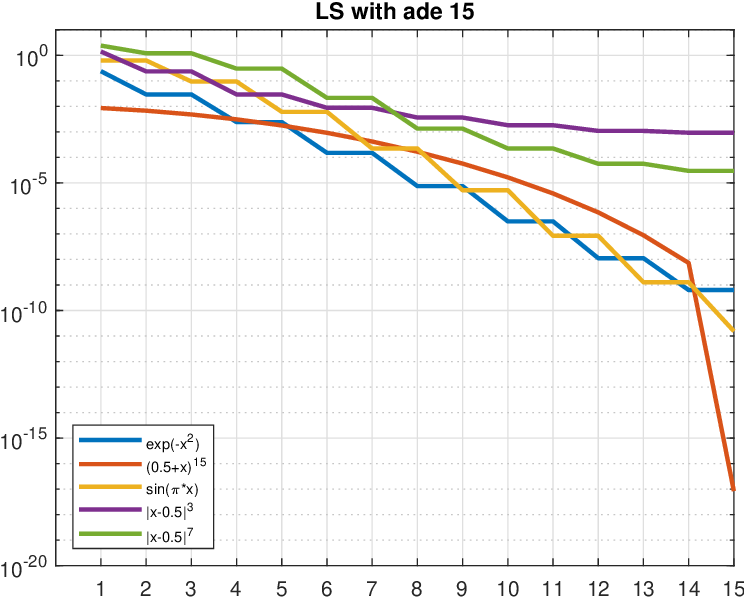}
\caption{{\color{black}Numerical comparisons on the unit-interval $[-1,1]$.  Hyperinterpolation error at degrees $1,2,\dots,15$ by Clenshaw-Curtis rule with $ADE=15$ (left), by Gauss-Legendre rule with $ADE=30$ (center), Least-Squares error using Clenshaw-Curtis rule with $ADE=15$ (right).}}
\label{figure:MZ_interval_ftest}
\end{figure}

\begin{figure}[ht!]
\centering
\includegraphics[scale=0.42,valign=c]{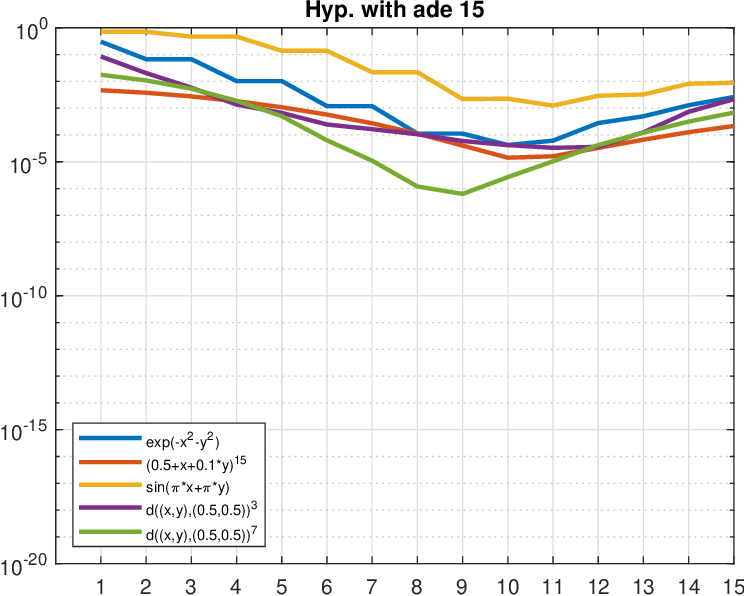}
\includegraphics[scale=0.42,valign=c]{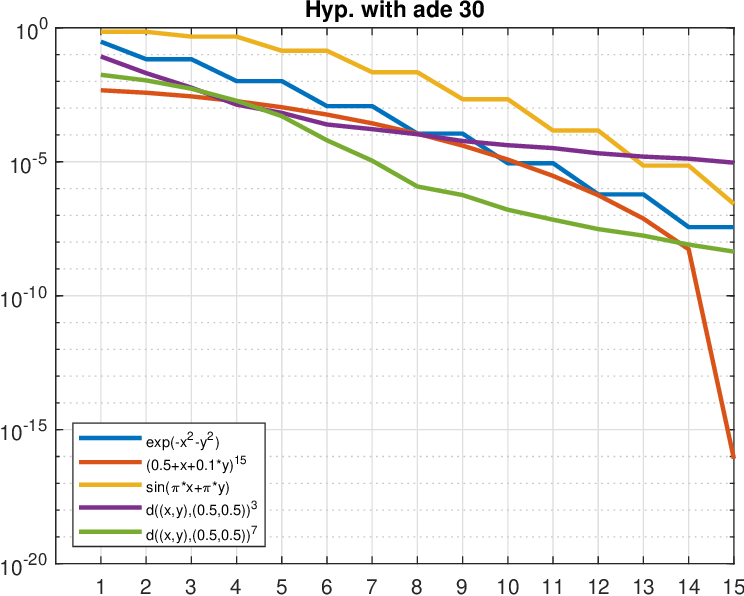}
\includegraphics[scale=0.42,valign=c]{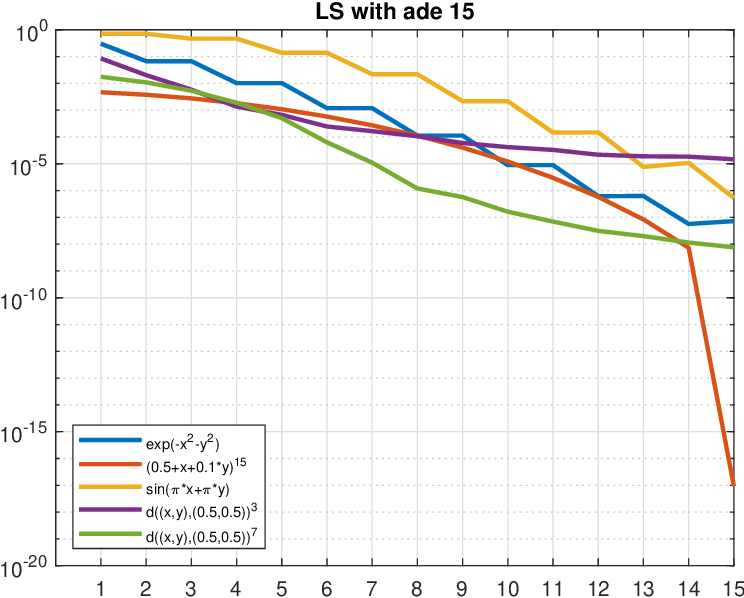}
\caption{{\color{black}Numerical comparisons on the unit-square $[-1,1]^2$. Hyperinterpolation error at degrees $1,2,\dots,15$ by Padua-points rule with $ADE=15$ (left), by Gauss-Legendre tensorial rule with $ADE=30$ (center), Least-Squares error using Padua-points rule with $ADE=15$ (right).}}
\label{figure:MZ_square_ftest}
\end{figure}

Finally we tested the approximation based on QMC rules on the unit-cube $[-1,1]^3$ of the functions
$$
f_1(x,y,z)=\exp(-(x^2+y^2+z^2)), \, f_2(x,y,z)=(0.5+x+0.1 y+0.4 z)^{15}, \, f_3(x,y,z)=\sin(\pi x + \pi y + \pi z),
$$
$$
f_4(x,y,z)=d((x,y,z),(0.5,0.5,0.5))^{3}=((x-0.5)^2+(y-0.5)^2+(z-0.5)^2)^{3/2}, 
$$
$$
\,f_5(x,y,z)=d((x,y,z),(0.5,0.5,0.5))^{7}=((x-0.5)^2+(y-0.5)^2+(z-0.5)^2)^{7/2}\;.
$$

The numerical results are reported in Figure {\ref{figure:MZ_square_ftest}, 
corresponding to QMC rule with $2^{15}=32768$ Halton points. One notices that QMC LS have a performance close to that of classical hyperinterpolation, while again QMC unfettered hyperinterpolation exhibits a  degradation increasing the polynomial degree $n$.

\begin{figure}[ht!]
\centering
\includegraphics[scale=0.42,valign=c]{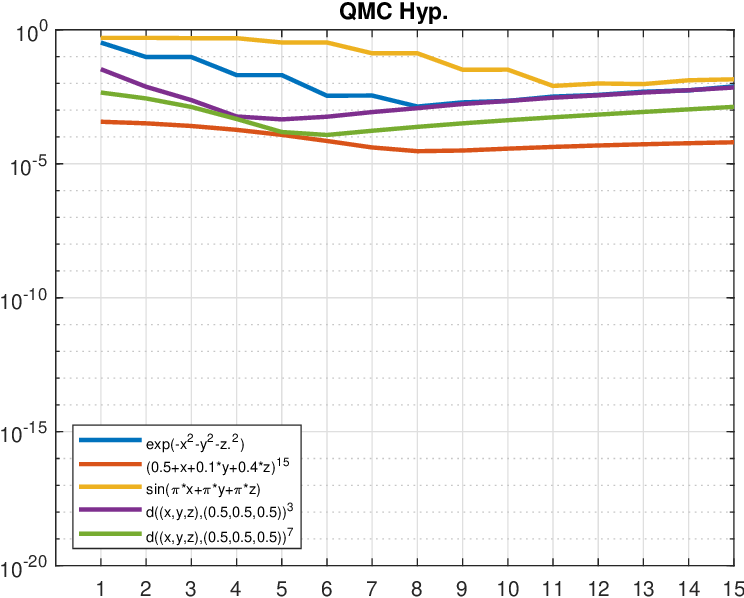} 
\includegraphics[scale=0.42,valign=c]{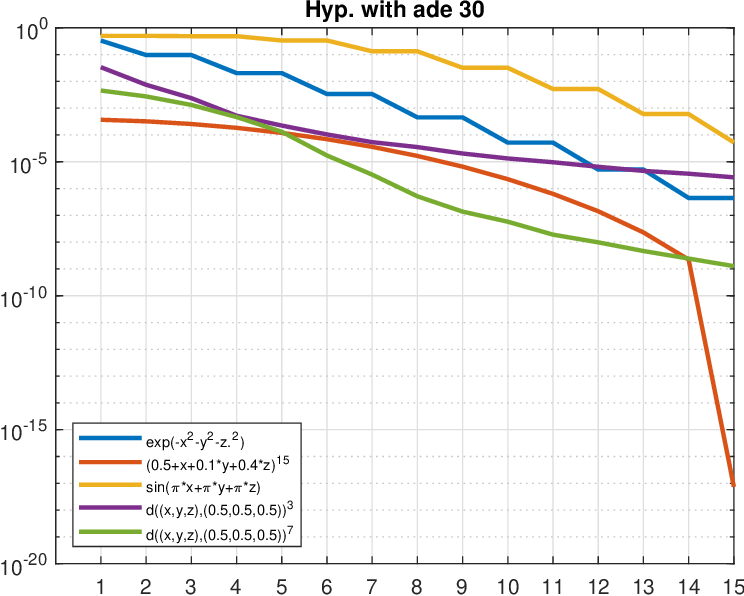} 
\includegraphics[scale=0.42,valign=c]{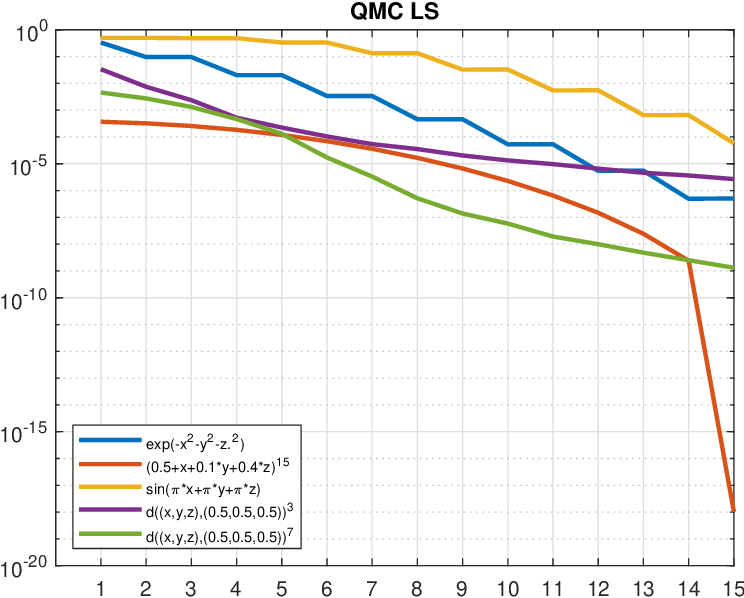}
\caption{{\color{black}Numerical comparisons on the unit-cube $[-1,1]^3$. QMC-hyperinterpolation errors by $2^{15}=32768$ Halton points (left), classical hyperinterpolation errors by Gauss-Legendre tensorial rule with $ADE=30$ (center), Least-Squares errors using QMC with $2^{15}$ Halton points (right).}}
\label{figure:MZ_cube_ftest}
\end{figure}
\vskip0.3cm
\noindent
{\bf Acknowledgements.} 
Work partially supported by the DOR funds of the University of Padova and by the INdAM-GNCS 2025 Project ``Polynomials, Splines and
Kernel Functions: from Numerical Approximation to Open-Source Software''. 
This research has been accomplished within the Community of Practice 
``Green Computing" of the Arqus European University Alliance, the RITA ``Research ITalian network on Approximation", and the SIMAI Activity Group ANA\&A. The work was partially supported by the National Natural Science Foundation of China (Project No. 12371099).

\end{document}